\documentclass[final,3p,times]{elsarticle}
\usepackage{graphicx}
\usepackage{amssymb}
\usepackage{amsmath}
\usepackage{subfig}
\usepackage{graphicx}
\usepackage{url}
\usepackage{algorithm}
\usepackage{algpseudocode}
\usepackage{multirow}
\usepackage{wasysym}
\usepackage{marvosym}
\usepackage{color, xcolor, colortbl}
\usepackage{array, booktabs}
\usepackage{xspace}
\usepackage{enumerate}
\usepackage[shortlabels]{enumitem}
\usepackage{listings} % Required package
\usepackage{bm}

\newcommand{\secref}[1]{Section \ref{#1}}
\newcommand{\fref}[1]{Fig.~\ref{#1}}

\renewcommand{\eqref}[1]{Eq.~(\ref{#1})}
\newcommand{\eref}[1]{(\ref{#1})}

\newcommand{\mat}[1]{\textrm{\textbf{#1}}}

\newcommand{\xten}[1]{$\times$ 10$^{\text{#1}}$}
\newcommand{\xtenm}[1]{\times \text{10}^{\text{#1}}}

\definecolor{light}{rgb}{0.8,0.8,0.8}
\definecolor{medium}{rgb}{0.6,0.6,0.6}
\definecolor{dark}{rgb}{0.4,0.4,0.4}
\definecolor{darkmed}{rgb}{0.3,0.3,0.3}
\definecolor{darkest}{rgb}{0.2,0.2,0.2}
\definecolor{Black}{rgb}{0,0,0}
\definecolor{White}{rgb}{1,1,1}
\definecolor{lightpurple}{rgb}{0.78823,0.709803,0.74509}
\definecolor{lightpurpletext}{rgb}{0.788235,0.5529411,0.658823}
\definecolor{skyblue}{rgb}{0.80392,0.866666,0.92941}
\definecolor{skybluetext}{rgb}{0.61568627,0.7647058,0.913725}
\definecolor{darkgreen}{rgb}{0.3137254,0.458823,0.18431}
\definecolor{foliagegreen}{rgb}{0.188,0.415,0.105}
\definecolor{steelbluegrey}{rgb}{0.1961,0.2353,0.2392}
\definecolor{highlightblue}{rgb}{0.4078,0.6431,0.85}
\definecolor{matlabblue}{rgb}{0,0.2705,0.85}
\definecolor{darkred}{rgb}{0.8,0.1725,0}
\definecolor{fireenginered}{rgb}{0.505,0.1411,0}
\definecolor{darkpurple}{rgb}{0.6431,0.3137,0.8509}
\definecolor{gaylordpurple}{rgb}{0.416,0.204,0.549}
\definecolor{deludedorange}{rgb}{0.7409,0.4392,0}
\definecolor{darksalmon}{rgb}{0.9137,0.411,0.706}
\begin{document}

\begin{frontmatter}



%~~~~~~~~~~~~~~~~~~~~~~~~~~~~~~~~~~~~
\title{Solving advection equations with reduction multigrids on GPUs\tnoteref{crown}}
\author[AMCG]{S. Dargaville}
\ead{dargaville.steven@gmail.com}
\tnotetext[crown]{UK Ministry of Defence © Crown owned copyright 2025/AWE}
\author[AWE]{R.P. Smedley-Stevenson}
\author[AMEC,AMCG]{P.N. Smith}
\author[AMCG]{C.C. Pain}
\address[AMCG]{Applied Modelling and Computation Group, Imperial College London, SW7 2AZ, UK}
\address[AWE]{AWE, Aldermaston, Reading, RG7 4PR, UK}
\address[AMEC]{ANSWERS Software Service, Amentum, Kings Point House, Queen Mother Square, Poundbury, DT1 3BW, UK}
%~~~~~~~~~~~~~~~~~~~~~~~~~~~~~~~~~~~~
\begin{abstract}
Methods for solving hyperbolic systems typically depend on unknown ordering (e.g., Gauss-Seidel, or sweep/wavefront/marching methods) to achieve good convergence. For many discretisations, mesh types or decompositions these methods do not scale well in parallel. In this work we demonstrate that the combination of AIRG (a reduction multigrid which uses GMRES polynomials) and PMISR DDC (a CF splitting algorithm which gives diagonally dominant submatrices) can be used to solve linear advection equations in parallel on GPUs with good weak scaling. We find that GMRES polynomials are well suited to GPUs when applied matrix-free, either as smoothers (at low order) or as an approximate coarse grid solver (at high order). To improve the parallel performance we automatically truncate the multigrid hierarchy given the quality of the polynomials as coarse grid solvers. Solving time-independent advection equations in 2D on structured grids, we find 66--101\% weak scaling efficiency in the solve and 47--63\% in the setup with AIRG, across the majority of Lumi-G, a pre-exascale GPU machine.
\end{abstract}
%~~~~~~~~~~~~~~~~~~~~~~~~~~~~~~~~~~~~
%~~~~~~~~~~~~~~~~~~~~~~~~~~~~~~~~~~~~
\begin{keyword}
Asymmetric multigrid \sep Parallel \sep GPU \sep AIRG \sep GMRES polynomials \sep PMISR DDC
\end{keyword}
%~~~~~~~~~~~~~~~~~~~~~~~~~~~~~~~~~~~~

\end{frontmatter}
%~~~~~~~~~~~~~~~~~~~~~~~~~~~~~~~~~~~~
%~~~~~~~~~~~~~~~~~~~~~~~~~~~~~~~~~~~~
%~~~~~~~~~~~~~~~~~~~~~~~~~~~~~~~~~~~~
%-----------------------------
\section{Introduction}
\label{sec:Introduction}
In this work we focus on solving the time-independent, constant direction advection equation with zero source in two dimensions, 
\begin{equation}
\mat{v} \cdot \nabla u = 0 \quad \textrm{in} \quad \Omega = [0, L_x] \times [0, L_y],
\label{eq:advection}
\end{equation}
where $u$ is the scalar quantity being advected, $\mat{v}$ is the constant direction velocity and $L_x$ and $L_y$ are the lengths of the domain. We apply an inflow boundary condition
\begin{equation}
u = g \quad \textrm{on} \quad \Gamma_{\textrm{in}}=\{\mat{x} \in \partial \Omega: \mat{v} \cdot \mat{n}(\mat{x})<0\},
\label{eq:advection_bc}
\end{equation}
where $g$ is a constant and $\mat{n}(\mat{x})$ is the boundary normal at $\mat{x}$. This is a canonical, linear PDE which underpins many linear and non-linear hyperbolic problems. 

The solution of \eref{eq:advection} is challenging to obtain given the requirement for a stable discretisation and the difficulties in solving the resulting asymmetric, non-normal, linear system which lacks strict diagonal dominance and whose condition number increases (like $\frac{1}{h}$) with increasing resolution. We focus on the time-independent equation as it is more difficult to solve than the time-dependent given that implicit time-stepping schemes typically improve the diagonal dominance, conditioning or symmetry (depending on the discretisation used). Unfortunately there is only a small subset of solver technologies, discretisations, mesh types and parallel decompositions which exhibit good weak scaling for advection equations with distributed memory parallelism (e.g., MPI). 
%-----------------------------
\subsection{Lower-triangular sparsity}
\label{sec:lower_triangular_sparsity}
%-----------------------------
Upwind Finite Difference (FD) and Discontinuous Galerkin (DG) Finite Element Methods (FEMs) can produce lower triangular structure in the discretised linear system and this property is often heavily exploited. The reason for this is the algorithmic efficiency of Gauss-Seidel methods. Information propagates across the mesh in the direction of the velocity (i.e., following the characteristics) and hence a single iteration of a Gauss-Seidel method following this ordering produces an exact solution. Algorithms following this principle are known by several names, including sweeps, wavefront or marching methods. In this paper we use the term ``sweep'' to refer to this algorithm, in contrast to a more general Gauss-Seidel iteration where an exact solution is not obtained in one iteration. (e.g., because a different ordering is used or there is a lack of lower triangular sparsity).  

Sweeps are an inherently sequential algorithm and can become challenging to perform in parallel as unknowns must wait for the wavefront to reach them before they can begin work (i.e., a scheduling problem). Some amount of lagging or asynchronous communication can be introduced in an attempt to obtain more parallelism \cite{Yavuz1992, Nowak1999, Rosa2010}, but this does not result in scalable methods, e.g., the number of iterations of a block-Jacobi method with a local sweep grows arbitrarily large when weak scaling. These methods can be used as preconditioners to outer iterations (e.g., GMRES) but this still does not give scalability. 
% ~~~~~~~~~~~~~
\begin{figure}[thp]
\centering
\subfloat[][Volumetric]{\label{fig:volumetric}\includegraphics[width =0.2\textwidth]{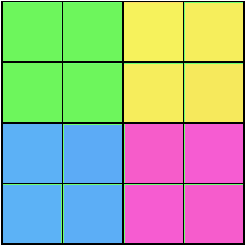}} \quad
\subfloat[][KBA]{\label{fig:pencil}\includegraphics[width =0.2\textwidth]{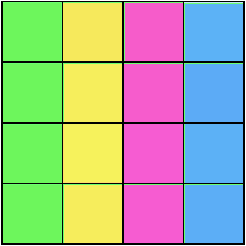}} \quad
\subfloat[][Overloaded]{\label{fig:overloaded}\includegraphics[width =0.2\textwidth]{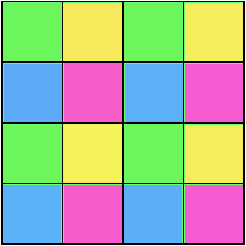}} \quad
\caption{Parallel decompositions for structured 2D mesh onto 4 MPI ranks (colours).}
\label{fig:decompositions}
\end{figure}
% ~~~~~~~~~~~~~

For structured meshes with upwind FD and DG FEM discretisations, the number of non-zeros per row is constant (giving a fixed communication pattern) and the choice of parallel decomposition (or partitioning) is key to attaining parallelism with a sweep. \fref{fig:decompositions} shows three common choices. The volumetric decomposition shown in \fref{fig:volumetric} minimises the edge cut of the adjacency graph (and hence minimises communication volume) and is often used in multi-physics. In advection problems however, it is a poor choice with a sweep as there is no point at which all ranks have work to do (e.g., consider $\mat{v}=(1,1)$ where a sweep will progress from bottom left to top right). The other two decompositions in \fref{fig:decompositions} provide more parallelism, at the cost of increased (though still scalable) communication volume. 

The KBA decomposition (or ``pencil'' decomposition) in \fref{fig:pencil} uses long thin columns (or pipes) \cite{Baker1998, Adams2013}. The parallel efficiency that can be achieved is a function of the resolution applied in each dimension. There must be enough pencils to distribute across all of the MPI ranks, with the length of the pencil determining the relative impact of the start-up phase (or pipe-fill) where not all ranks have work to do. In 2D with equal resolution applied in both dimensions, the best parallel efficiency that can be achieved with this decomposition is 50\%; there is only one ``stage'' where all ranks are active (when the wavefront reaches the diagonal). On a structured mesh the impact of the pipe-fill is greater in 2D than in 3D \cite{Pautz2017}, making it more challenging to achieve good parallel efficiency in 2D. 

The ``overloaded'' decomposition \cite{Pautz2017} in \fref{fig:overloaded} is nested and non-contiguous, in an attempt to give ranks work throughout the sweep process. It decomposes the mesh in two steps, the first typically into fewer partitions than ranks, with the second taking each partition generated by the first and decomposing into as many partitions as ranks. \cite{Pautz2017} shows similar scaling on structured meshes in 3D with an overloaded decomposition as with a KBA. The non-contiguous partitioning would however make saturating local, shared memory compute difficult (e.g., GPUs).

In Boltzmann applications such as radiation transport, spectral-wave modelling, Lattice-Boltzmann and kinetic theory, the discretised streaming operator gives a linear system with block structure, where each block is an advection equation for a different direction. This can be exploited to hide the start-up phase by not waiting for each single direction sweep to finish before starting another (``pipelining''). Given this, KBA decompositions \cite{Baker1998} and related variants achieve near ideal weak scaling on both CPU and GPU based machines \cite{Baker2012a, Adams2013, Deakin2016, Kunen2019}. This is also the case for semi-structured meshes, where the local component of the mesh may be unstructured, but the parallel decomposition still requires structure \cite{Adams2020}. The downside of pipelining is that many directions are required with which to pipeline; in Boltzmann applications this requirement is not difficult to meet as there are often tens to hundreds of directions available per octant given the combination of the angular/energy dimensions.

If we consider general unstructured meshes, achieving good parallel scaling with sweeps is very difficult \cite{Pautz2002, Plimpton2005, Colomer2013} and calculating an optimal sweep schedule is NP complete. Upwind DG FEMs are typically used, and the connectivity of unstructured meshes means the number of non-zeros per row is not constant. Cycles can also appear in the mesh (i.e., loops in the directed graph), which destroy the lower-triangular structure \cite{Vermaak2021}. We can also no longer create KBA decompositions, though there have been attempts to impose the required structure on unstructured meshes (e.g., see Figure 9 \& 10 in \cite{Vermaak2021}, or \cite{Plimpton2005}). Both the volumetric and overloaded decompositions can be performed on an unstructured mesh (e.g., with a tool like ParMETIS \cite{karypis_metis-unstructured_1995}), though only the overloaded retains much distributed parallelism. 
%-----------------------------
\subsection{Lacking lower-triangular sparsity}
\label{sec:lacking_sparsity}
%-----------------------------
The specific combination of discretisation (upwind FD or DG FEM), mesh type (structured), solver technology (sweeps) and application (Boltzmann with many directions) described above is perhaps the only one commonly used on large parallel machines. Other combinations are possible, but typically do not scale. For example, stabilised Continuous Galerkin (CG) FEMs applied to advection equations do not produce lower-triangular structure and hence sweeps cannot be used. A typical way to solve the resulting linear system would be to apply an LU factorisation or ILU(0) preconditioner with an outer GMRES, but LU/ILU factorisations require backward/forward substitution which are Gauss-Seidel/sweep methods with the same limitations described above. Multigrid/multilevel methods are often used as solvers/preconditioners in elliptic problems but have traditionally struggled with asymmetric linear systems and often use ILU or Gauss-Seidel smoothers for robustness. 
%-----------------------------
\subsection{Asymmetric multigrids}
\label{sec:asymmetric_multigrids}
%-----------------------------
Over the last several years however, theory has been developed for multigrid in asymmetric linear systems \cite{Southworth2017, Manteuffel2019, Manteuffel2019a} and new algebraic multigrid methods have emerged, known as reduction multigrids. Variants of these methods have shown excellent convergence in common asymmetric equations with a range of discretisations, such as advection-diffusion, anisotropic diffusion, recirculating flow, and streaming operators from Boltzmann applications \cite{Southworth2017, Manteuffel2019, Manteuffel2019a, Sivas2021, Ali2024, Zaman2024, Dargaville2024a, Dargaville2024}. Importantly these methods do not rely on any specific mesh type, smoother, unknown ordering, lower-triangular or block structure for good convergence. 

The parallel scaling of reduction multigrids has been investigated on CPU architectures in both \cite{Hanophy2020} and our previous work \cite{Dargaville2025}, specifically with streaming operators from Boltzmann applications. Our goal in \cite{Dargaville2025} was to show that the combination of a reduction multigrid, AIRG (approximate ideal restriction with GMRES polynomials), and a CF splitting algorithm, PMISR DDC, scaled well on unstructured grids decomposed with a volumetric decomposition, using a stabilised CG FEM discretisation that lacked lower triangular structure. Unfortunately our ability to generate/decompose fully unstructured meshes limited the weak scaling results we could show. Given this, we attempted to expose communication bottlenecks by using very few DOFs/core and only used 4 directions to try and demonstrate that parallelism comes from the spatial dimension, rather than pipelining many directions like in sweeps. We showed 81\% weak scaling efficiency in the solve from 256 cores to 8196 cores with only 8.8k DOFs/core.
%-----------------------------
\subsection{Aims}
\label{sec:aim}
%-----------------------------
One of the aims of this work is to follow on from \cite{Dargaville2025} and demonstrate that excellent parallel performance is possible in advection-type problems without sweeps at larger scales. We also wish to show reduction multigrids are well suited to GPUs. We originally designed AIRG and the PMISR DDC CF splitting to be suitable for use on different hardware and in this work we detail several changes (compared to \cite{Dargaville2025}) that improve performance on GPUs, namely:
\begin{enumerate}
\item We use low-order GMRES polynomials applied matrix-free as our F-point smoothers
\item We use high-order GMRES polynomials in the Newton basis applied matrix-free as our coarse grid solver
\item We automatically truncate the multigrid hierarchy given the efficacy of the coarse-grid solver
\end{enumerate}

We show results from running on Lumi-G, a pre-exascale GPU machine. Ideally we would like to run with both structured and fully unstructured meshes for a range of discretisations in 2D and 3D. This is difficult however, for example in 2D we would be required to generate/decompose an unstructured mesh with close to a trillion elements to scale to the entirety of Lumi-G. Although they are difficult to generate at scale, reduction multigrids in advection problems perform best on unstructured grids, due to the lack of consistent alignment between the velocity direction and mesh faces. Structured meshes, while easy to generate/decompose at scale, represent a worst case for advection problems solved with reduction multigrids. They are typically harder to solve than with unstructured meshes, requiring slower coarsenings and hence more levels (and more memory), with convergence dependent on the velocity direction (e.g., see Figure 4 in \cite{Manteuffel2019}). 

Given these factors, we challenge our methods and use a single direction and structured quadrilateral meshes in 2D in this paper. We use a simple upwind FD discretisation and decompose our mesh with a volumetric decomposition. We emphasize that we do not take advantage of the mesh structure or sparsity at any point (e.g., we do not order the unknowns to explicitly give lower-triangular structure in the matrix), our methods are entirely algebraic. We show that a GPU port of our methods gives excellent weak scaling in the solve on Lumi-G. The code used in this work (AIRG, the PMISR DDC CF splitting, along with other reduction multigrids such as nAIR and $\ell$AIR \cite{Southworth2017, Manteuffel2019}) is available in the open-source \textsc{PFLARE}\footnote{\url{https://github.com/PFLAREProject/PFLARE}} library, which can be installed directly through the PETSc configure (\texttt{--download-pflare}). 

We are not aware of any method that shows the level of parallel scaling achieved in this work for advection (or related) problems. Given this and the wider applicability of these methods to other asymmetric problems, we believe reduction multigrids bring parallel capabilities to many fields which were previously lacking them.
%-----------------------------
\section{Reduction multigrids}
\label{sec:reduction}
%-----------------------------
We start with a coarse/fine (CF) splitting of the unknowns, designed to produce a fine-fine block that allows a good quality sparse inverse, while coarsening as fast as possible. If we write a linear system $\mat{A}\mat{x}=\mat{b}$ of size $n \times n$ as a block system given by the coarse and fine points $n_\textrm{C} + n_\textrm{F} = n$ we have
\begin{equation}
\begin{bmatrix}
\mat{A}_\textrm{ff} & \mat{A}_\textrm{fc} \\
\mat{A}_\textrm{cf} & \mat{A}_\textrm{cc}
\end{bmatrix}
\begin{bmatrix}
\mat{x}_\textrm{f} \\
\mat{x}_\textrm{c}
\end{bmatrix} = 
\begin{bmatrix}
\mat{b}_\textrm{f} \\
\mat{b}_\textrm{c}
\end{bmatrix}.
\label{eq:air_two}
\end{equation}
Ideal prolongation and restriction operators can be written as
\begin{equation}
\mat{P} = 
\begin{bmatrix}
\mat{W} \\
\mat{I}
\end{bmatrix}, \quad
\mat{R} = 
\begin{bmatrix}
\mat{Z} & \mat{I}
\end{bmatrix},
\label{eq:prolong}
\end{equation}
and approximate ideal restriction (AIR) or approximate ideal prolongation (AIP) operators can be formed by seeking approximate solutions to 
\begin{equation}
\quad \mat{Z} \mat{A}_\textrm{ff} =-\mat{A}_\textrm{cf} \quad \textrm{and} \quad \mat{A}_\textrm{ff} \mat{W} = -\mat{A}_\textrm{fc}.
\label{eq:z}
\end{equation}
The coarse-grid matrix can be calculated by the matrix triple-product $\mat{A}_\textrm{coarse}=\mat{R}\mat{A}\mat{P}$ and this process is then repeated to generate a multigrid hierarchy. Different reduction multigrids, such as nAIR, $\ell$AIR and AIRG are characterised by how they solve \eref{eq:z} and/or the F-point smoothing they apply.

In this work we compute an approximate ideal restrictor by constructing an approximation $\hat{\mat{A}}_\textrm{ff}^{-1} \approx \mat{A}_\textrm{ff}^{-1}$ and then performing a matrix-matrix product \eref{eq:z}, i.e., $\mat{Z} =-\mat{A}_\textrm{cf} \hat{\mat{A}}_\textrm{ff}^{-1}$. As in \cite{Manteuffel2019, Dargaville2025} we use a classical one-point prolongator rather than an ideal approximate prolongator; \cite{Manteuffel2019} note that either an ideal restrictor or ideal prolongator is required, not both.
 
In our previous work \cite{Dargaville2024a, Dargaville2024, Dargaville2025} we used the same $\hat{\mat{A}}_\textrm{ff}^{-1}$ to perform only up F-point smoothing on each level (i.e., there is no down smoothing or smoothing on C points). In this work we also perform only up F-point smoothing, however we use a different (related) approximation, $\bar{\mat{A}}_\textrm{ff}^{-1} \approx \mat{A}_\textrm{ff}^{-1}$. Hence, on each level if we are performing F-point smoothing on $\mat{A}\mat{e}=\mat{r}$, where $\mat{r}$ is the residual computed after the coarse-grid correction, then
\begin{equation}
\mat{e}_\textrm{f}^{n+1} = \mat{e}_\textrm{f}^{n} + \bar{\mat{A}}_\textrm{ff}^{-1}(\mat{r}_\textrm{f} - \mat{A}_\textrm{fc} \mat{e}_\textrm{c}^{n} - \mat{A}_\textrm{ff} \mat{e}_\textrm{f}^{n}). 
\label{eq:f_point}
\end{equation}
If multiple F-point smooths are performed in a row the value of $\mat{A}_\textrm{fc} \mat{e}_\textrm{c}^{n}$ can be stored as it does not change. We also add relative drop tolerances to $\mat{A}_\textrm{coarse}$ (we lump the entries dropped from $\mat{A}_\textrm{coarse}$ to the diagonal) and $\mat{R}$ once they have been assembled to help keep the complexities and hence memory use down. We can also see from \eref{eq:f_point} that we don't have to store $\mat{A}_\textrm{coarse}$ during the solve if we are only performing F-point smooths. After we have built the coarse matrix on each level, we only need to keep $\mat{A}_\textrm{ff}$ and $\mat{A}_\textrm{fc}$. This further helps reduce the memory required. This is important for reduction multigrids, as compared to elliptic multigrid methods, the coarsening can be slow and hence there are many C points. 

If using an exact $\mat{A}_\textrm{ff}^{-1}$ instead of sparse approximations, reduction multigrids can become exact solvers. This is because the ideal restrictor ensures that after coarse-grid correction, the C-point error is zero. This leaves error on only F points, which is reduced to zero with an exact F-point smooth. Given this, it is often worth thinking of smoothers for reduction multigrid methods as solvers, in that we want all the error modes to go to zero, rather than just a subset (e.g., high frequency modes as in many elliptic multigrids). Outside of a few specific cases, this limit is not practical given that $\mat{A}_\textrm{ff}^{-1}$ is often dense. But unlike in many multigrid methods, there is at least a straightforward route to improving convergence in asymmetric problems with reduction multigrids, by improving the approximation of $\mat{A}_\textrm{ff}^{-1}$. 

One such special case is when $\mat{A}$ is tridiagonal (e.g., 1D diffusion/advection). Reduction multigrids in this case are equivalent to cyclic reduction methods \cite{Heller1976, Golub1992, Gander1998, Bini2009}, which perform well on GPUs and in parallel \cite{Balogh2022, Song2022, Tolmachev2025}. A red-black CF splitting can be performed on every level, giving two independent sets and hence $\mat{A}_\textrm{ff}$ and $\mat{A}_\textrm{cc}$ are diagonal and can be exactly inverted. Recursive hierarchies can then be built either for one of the sets of F or C points, or both to produce an exact solver. 

Another special case is when the CF splitting gives a maximal independent set in the adjancency graph of $\mat{A}$, by specifying a strong threshold of 0 in the CF splitting. This gives an independent set for the F points, but not the C points. This gives diagonal $\mat{A}_\textrm{ff}$ that can be exactly inverted, but also results in a very slow coarsening with many levels. This can impact the parallel performance and the memory use, although not as much as might be expected \cite{Dargaville2025}. 

In general however, reduction multigrids are designed to work without any assumptions about structure. The only requirement for good performance is to produce a CF splitting that gives a large $\mat{A}_\textrm{ff}$ whose inverse can be sparsely approximated. If the approximation to $\mat{A}_\textrm{ff}^{-1}$ is poor, then the C-point error after coarse grid correction will not be close to zero. C-point smoothing can be added on each level to help reduce this error, but this means storing $\mat{A}_\textrm{cc}$ and $\mat{A}_\textrm{cf}$ (i.e., on each level the entirety of $\mat{A}_\textrm{coarse}$ must be kept to perform both F and C point smoothing). Balancing the different choices for computing the CF splitting and $\mat{Z}$, performing F and/or C point smoothing and the resulting memory requirements can be difficult, particularly in parallel and on differing hardware (e.g., CPUs and GPUs). In the following sections we describe the choices we make and their suitability for GPUs. 
%-----------------------------
\subsection{CF splitting}
\label{sec:cf_splitting}
%-----------------------------
In the previous section we noted that reduction multigrids rely on the coarse/fine splitting producing an $\mat{A}_\textrm{ff}$ block that is ``easy'' to sparsely approximate. Specifically we aim to find a large $\mat{A}_\textrm{ff}$ submatrix that is diagonally dominant. In \cite{Dargaville2025} we introduced a two-pass CF splitting algorithm designed to produce diagonally dominant $\mat{A}_\textrm{ff}$, called PMISR DDC (Parallel Modified Independent Set for Reduction multigrids and Diagonal Dominance Cleanup). \cite{Dargaville2025} showed that using the PMISR DDC CF splitting with AIRG produced the best convergence when tested against common CF splittings such as PMIS, HMIS and Falgout-CLJP, while also being the most performant. 

If we denote $\mat{S}$ as an adjacency matrix which defines the strong connections for each unknown, PMISR computes an independent set in $\mat{S} + \mat{S}^\textrm{T}$ with a Luby-type algorithm \cite{Luby1985}, ensuring for an asymmetric linear system that there are no large off-diagonal entries in $\mat{A}_\textrm{ff}$. In \cite{Dargaville2025} we set a fixed number of Luby steps per level in parallel ($n_\textrm{loops}^\textrm{max}=3$) to limit the number of synchronisation points. In this work we allow as many as needed for robustness across a range of problems, but note that for any given problem, the CF splitting time can be reduced by fixing $n_\textrm{loops}^\textrm{max}$. 

Many small entries in $\mat{A}_\textrm{ff}$ could still cause a loss of diagonal dominance, so the second pass, DDC, converts F points to C points based on a diagonal dominance ratio. In \cite{Dargaville2025} we performed one DDC pass which converted a fixed percentage (10\%) of the least diagonally dominant local F points to C points. In this work we allow multiple DDC passes (given they are cheap) and convert a smaller fixed percentage of local F points per pass (1\%), as this results in a slightly smaller grid complexity. 

Given we only use structured, quadrilateral grids with an upwind FD stencil in this paper, we could generate a red-black CF splitting on the top grid, regardless of the velocity direction. We chose not to do this so that our method stays entirely algebraic, but note that this would result in a faster coarsening and hence better performance. 

One of our goals in \cite{Dargaville2025} was to produce a CF splitting algorithm that was compatible with different hardware. Both the PMISR and DDC algorithms are suitable for use on GPUs. The PMISR algorithm is equivalent to a slightly modified PMIS \cite{Sterck2006} and hence very parallel. The DDC algorithm computes the diagonal dominance ratio for each row, which is embarrassingly parallel and then computes the target ratio which would result in the requested fixed percentage of F points be converted to C points. This can be calculated by binning the diagonal dominance ratios and choosing the closest bin boundaries. We use 1000 bins and note that as long as there is a range of diagonal dominance ratios there is very little contention in this algorithm. On a structured grid, all rows have the same diagonal dominance ratio on the top grid (excluding the boundaries) and hence we would consider the examples in this paper to be a worst case. Despite this, the cost of DDC passes is small and hence we do not consider the parallelism of the DDC algorithm any further. We also use work overlapped asynchronous communication where possible in both the PMISR and DDC passes. Given we can produce a digonally dominant $\mat{A}_\textrm{ff}$ in parallel on GPUs, the next section discusses how we approximate $\mat{A}_\textrm{ff}^{-1}$. 
%-----------------------------
\subsection{AIRG}
\label{sec:airg}
%-----------------------------
Approximate ideal restriction with GMRES polynomials (AIRG) uses fixed-order GMRES polynomials \cite{saad_gmres:_1986, Nachtigal1992, Loe2022} to compute approximations to $\mat{A}_\textrm{ff}^{-1}$ \cite{Dargaville2024a}. In the GMRES algorithm, if we consider the F-point linear system $\mat{A}_\textrm{ff}\mat{x}_\textrm{f}=\mat{b}_\textrm{f}$, the solution at step $m$ can be written as $\mat{x}^m_\textrm{f} = q_{m-1}(\mat{A}_\textrm{ff}) \mat{b}_\textrm{f}$, where $q_{m-1}(\mat{A}_\textrm{ff})$ is a matrix polynomial of degree $m-1$ known as the GMRES polynomial and is given by
\begin{equation}
q_{m-1}(\mat{A}_\textrm{ff}) = \alpha_0 + \alpha_1 \mat{A}_\textrm{ff} + \alpha_2 \mat{A}_\textrm{ff}^2 + \ldots + \alpha_{m-1}\mat{A}_\textrm{ff}^{m-1} \approx \mat{A}_\textrm{ff}^{-1},
\label{eq:gmres_poly}
\end{equation}
where $\alpha$ are the polynomial coefficients which minimise the 2-norm of the residual at step $m$ and can be output by a (slightly modified) GMRES method. Typically GMRES methods do not explicitly compute the polynomial coefficients, as the polynomial can be applied during the orthogonalisation in a numerically stable way. 

Rather than use GMRES directly (which would require dot products during the solve), we generate a random rhs and calculate the polynomial coefficients once per level in the setup. This stationary polynomial is then used as our approximation $\bar{\mat{A}}_\textrm{ff}^{-1} = q_{m-1}(\mat{A}_\textrm{ff})$, as it can be applied with only Sparse Matrix-Vector (SpMV) products (i.e., no dot products). Stably computing the polynomial coefficients can be difficult, but we only require low-order polynomials (small $m$) for good approximations, given the $\mat{A}_\textrm{ff}$ on each level are diagonally dominant. We use an Arnoldi method to compute them for robustness \cite{Nachtigal1992}. Previously we used a communication-avoiding method that used the power basis \cite{Dargaville2025}, but in this work we found that even with diagonally dominant $\mat{A}_\textrm{ff}$ and low-order polynomials, this method was less stable than the Arnoldi at the largest scale, without a significant difference in relative cost. 

These polynomials give good approximations with symmetric or asymmetric $\mat{A}_\textrm{ff}$ while also being agnostic to unknown ordering, block and/or triangular structure. We noted in \secref{sec:reduction} we use two related approximations to $\mat{A}_\textrm{ff}^{-1}$, namely $\bar{\mat{A}}_\textrm{ff}^{-1}$ in the F-point smoothing and $\hat{\mat{A}}_\textrm{ff}^{-1}$ to build the approximate ideal restrictor. We use $\bar{\mat{A}}_\textrm{ff}^{-1} = q_{m-1}(\mat{A}_\textrm{ff})$ as our F-point smoother, as we can apply the polynomials during the solve matrix-free. This only requires the storage of the polynomial coefficients, which are computed once in the setup. We could also use a matrix-powers kernel to do this with fewer communication steps, but this would take extra memory so we chose not to in this work.

For the other approximation we use $\hat{\mat{A}}_\textrm{ff}^{-1} \approx q_{m-1}(\mat{A}_\textrm{ff})$. We use the same polynomial coefficients, but we assemble an approximation of $q_{m-1}$ with controlled sparsity so we can cheaply compute the matrix-matrix product $\mat{Z} =-\mat{A}_\textrm{cf} \hat{\mat{A}}_\textrm{ff}^{-1}$. As described in \cite{Dargaville2024a, Dargaville2025}, we enforce the sparsity of a given matrix-power on $\hat{\mat{A}}_\textrm{ff}^{-1}$, namely the sparsity of $\mat{A}_\textrm{ff}$, or
\begin{equation}
\hat{\mat{A}}_\textrm{ff}^{-1} = \alpha_0 + \alpha_1 \mat{A}_\textrm{ff} + \alpha_2 \tilde{\mat{A}}_\textrm{ff}^2 + \ldots + \alpha_{m-1} \tilde{\mat{A}}_\textrm{ff}^{m-1},
\label{eq:fixed_sparsity}
\end{equation}
where each subsequent matrix power greater than one has entries outside the sparsity of $\mat{A}_\textrm{ff}$ dropped. This allows us to assemble our approximation $\hat{\mat{A}}_\textrm{ff}^{-1}$ very cheaply with less communication, while still retaining good performance in these problems. This kernel is well suited to GPUs, as we are performing a matrix-matrix product, but one where subsequent powers share the same sparsity.

Previously \cite{Dargaville2024a, Dargaville2025} we used the same assembled approximation $\hat{\mat{A}}_\textrm{ff}^{-1}$ for both the F-point smoothing and the restrictor. The choice to separate the approximations in this work is due to our use of GPUs. $\bar{\mat{A}}_\textrm{ff}^{-1}$ is a better approximation of $\mat{A}_\textrm{ff}^{-1}$, requires no extra memory but needs $m$ matrix-vector products to apply. On the other hand $\hat{\mat{A}}_\textrm{ff}^{-1}$ is a poorer approximation of $\mat{A}_\textrm{ff}^{-1}$, requires storing an assembled matrix but only needs one matrix-vector product to apply. 

For good performance on GPUs we need as many DOFs per rank as possible, hence we would like to reduce the memory used given the slow coarsening required in asymmetric linear systems. Applying the polynomials matrix-free increases the amount of work as noted, but decreases memory use. Furthermore, the SpMV products used to apply the polynomial are memory-bandwidth bound, not FLOP bound and hence GPUs are suited to an approximation which repeatedly applies the same matrix. We see this reflected in the solve time. For example, we use $6^\textrm{th}$ order polynomials matrix free as our F-point smoothers in \secref{sec:Results}. In the 2 node problem shown in \fref{fig:strong_scaling}, this results in 7 iterations and a per V-cycle time of 0.11s. This requires six SpMVs on each level in our F-point smoothing to apply $\bar{\mat{A}}_\textrm{ff}^{-1}$. If instead we use the assembled approximation $\hat{\mat{A}}_\textrm{ff}^{-1}$ in our F-point smoothing, this only requires one SpMV. We find the iteration count increase to 8, with the per V-cycle time however only slightly decreasing to 0.090s. So although we perform far more work applying our polynomials matrix-free, this work is well suited to GPUs and we benefit from improved convergence. As such we always apply our polynomials matrix-free when F-point smoothing in the results shown below.
%-----------------------------
\subsection{Coarse grid solver}
\label{sec:coarse grid}
%-----------------------------
In \cite{Dargaville2025} we applied several iterations of an assembled, low-order GMRES polynomial as the coarse grid solver and truncated the hierarchy to try and avoid communication bottlenecks on the coarse grids. Given the GPU hardware, in this paper we instead use high-order GMRES polynomials applied matrix-free as the coarse grid solver. For stability, we compute the harmonic Ritz values (which are the roots of the GMRES residual polynomial) and apply the polynomial as in \cite{Loe2022}. Applying the GMRES polynomial as a Newton polynomial with the Harmonic Ritz values as roots is equivalent to using a Newton basis, but without the need to explicitly store a basis (we can just store the roots). We take care to use a rank-revealing factorisation when computing the harmonic Ritz values for stability. For example if we were trying to compute the harmonic Ritz values for a scaled identity (like if we had performed a red-black CF splitting on a structured grid, giving diagonal $\mat{A}_\textrm{ff}$) then a low-order polynomial would be an exact inverse and attempting to compute higher order terms can be unstable. We find this combined with the added roots described by \cite{Loe2022} (and being careful about overflow in the product of factors) can allow us to apply GMRES polynomials stably up to very high order ($1000^\textrm{th}$ order in some problems).

There are several reasons why high-order GMRES polynomials applied matrix-free are well suited for use as coarse grid solvers, particularly on GPUs. They dampen all error modes, unlike some iterative methods which only smooth high-frequency errors. As mentioned in the previous section, their application is low-memory and only requires SpMVs (i.e., no dot products) which are suited to GPU hardware. Given that we can generate high-order polynomials, we can also use them to heavily truncate the multigrid hierarchy without affecting the convergence of the method. This is key to avoiding communication bottlenecks on the lower levels, where there is very little work in the setup/solve to hide communication. Truncating the hierarchy reduces the memory use, but also means the coarse grid matrix is larger and hence the amount of work in the solve goes up. Again however, GPUs are well suited to an approximation which trades extra work for less communication bottlenecks. This also avoids the problem of large relative kernel launch times on GPUs on small coarse grids, particularly if they have been redundantly distributed in an attempt to decrease the solve time (e.g., with PCREDUNDANT in PETSc). Multigrid methods often pin some of the lower levels to the CPUs (e.g., see GAMG in PETSc), but this has the cost of data transfers between the (separate) host and device on many architectures. Our approach instead keeps all of the levels on the GPU. 
%-----------------------------
\subsection{Automatic truncation}
\label{sec:truncation}
%-----------------------------
One of the practical difficulties in \cite{Dargaville2025} was choosing what level to truncate the hierarchy at, given an iterative coarse grid solver. As the cost of computing and applying the GMRES polynomial matrix-free can be relatively small, we now allow automatic truncation of the hierarchy. On each level during the setup, before computing the CF splitting, we generate a high-order GMRES polynomial to use as a tentative coarse grid solver, $q_{m_\textrm{coarse}-1}(\mat{A}_\textrm{coarse})$. We then test if this polynomial can solve $\mat{A}_\textrm{coarse} \mat{x} = \mat{b}$ to a given relative tolerance, $\alpha_\textrm{trunc}$. We set $\mat{b}$ to be a random vector (n.b., a different random vector than was used to generate the polynomial coefficients). If the tolerance $\alpha_\textrm{trunc}$ is reached, we truncate the hierarchy on that level and use the polynomial as the coarse grid solver, if not we discard the polynomial and continue coarsening. To reduce the cost, we can specify the level at which to start testing tentative coarse grid solvers, which is not difficult to estimate \textit{a priori} given we know the number of unknowns and an upper bound on the coarsening rate in advection problems (around two). In this work we find a very loose relative tolerance of $\alpha_\textrm{trunc}=1 \xtenm{-1}$ allows heavy truncation, while giving the same iteration count in the solve as using no truncation. We discuss the choice of the coarse grid solver polynomial order, $m_\textrm{coarse}$, in \secref{sec:coarse_order}. 
%~~~~~~~~~~~~~~~~~~~~~~~~~~~~~~~~~~~~
\section{Results}
\label{sec:Results}
% ~~~~~~~~~~~~~
To illustrate the performance of our methods, we solve \eref{eq:advection} subject to \eref{eq:advection_bc} on a structured 2D quadrilateral mesh, discretised with an upwind FD method. Unless otherwise noted, we use a square domain with $L_x = L_y = 1$, apply uniform resolution in the $x$ and $y$ directions, set the velocity to $\mat{v}=(\cos(\pi/4), \sin(\pi/4))$ and use a zero inflow boundary condition ($g=0$) on the left and bottom boundaries. We non-dimensionalise the problem with mesh resolution so we retain the same stencil values with refinement, giving the same matrix as the advection 2D gallery example in PyAMG \cite{pyamg2023} (although we do not eliminate the boundary DOFs). This non-dimensionalisation is equivalent to a diagonal scaling of our matrix and we note that diagonal scaling and/or block-element scaling can be important for performance, depending on the discretisation (e.g., see \cite{Manteuffel2019, Southworth2017}). 

We use our reduction multigrid as a solver, applied with an undamped Richardson to a relative tolerance of $1 \xtenm{-10}$ (measured with an unpreconditioned norm) and initial guess of one. We measure both the cycle and storage complexity of our multigrids. The cycle complexity is a count of the floating-point operations (FLOPs) required to perform one multigrid V-cycle, scaled by the number of FLOPs required to apply the top grid matrix. The storage complexity is defined in \cite{Dargaville2025} and reflects the amount of storage required in the solve: it is the sum of number of non-zeros in each matrix required during the solve, scaled by the number of non-zeros in the top grid matrix. Typically elliptic multigrids use the operator complexity as a measure of the storage, but this is not a good representation in a reduction multigrids given, for example, we only perform F-point smoothing in this work and hence do not need to store $\mat{A}_\textrm{coarse}$ on each level. We note that both metrics include the ``ghost'' entries on each rank, making them representative of the work/storage required in parallel. 

Typically multigrid methods are tested with a random initial condition (or rhs), but our methods are built on top of PETSc 3.23.1 \cite{petsc-web-page, mills2021} and calling \texttt{VecSetRandom} currently involves the host. Using a random vector would therefore trigger a host to device copy during our solve and we wish to avoid any copies when timing our methods. As such we use an initial guess of one but we verified our iteration count is the same with a random initial guess. As mentioned in \secref{sec:Introduction} the code used in this paper is available in the open-source \textsc{PFLARE} library. All the results generated below used PFLARE v1.22 and the default options of \texttt{tests/adv\_diff\_2d.c}. We trigger two solves with the command-line option \lstinline{-second_solve} and show the timing results from the second solve. This is because the discretisation of the advection equation in \texttt{tests/adv\_diff\_2d.c} is currently performed on the host and is copied to the device during the first solve. We use the Kokkos abstraction in PETSc to enable cross-platform GPU support and we wrote Kokkos implementations of AIRG, PMISR DDC and other required routines (e.g., for dropping entries according to a tolerance, performing F/C point smoothing, etc) in PFLARE to improve the setup/solve times. In \texttt{tests/adv\_diff\_2d.c} we use the command-line options \lstinline{-dm_mat_type aijkokkos -dm_vec_type kokkos} to specify the GPU implementation should be used. Unless otherwise noted, we use the following solver/preconditioner command-line options:
\begin{lstlisting}
-ksp_rtol 1e-10 -ksp_atol 1e-50 -ksp_type richardson -ksp_norm_type unpreconditioned -pc_type air -pc_air_a_lump -pc_air_a_drop 1e-6 -pc_air_strong_threshold 0.99 -pc_air_ddc_fraction 0.01 -pc_air_ddc_its 2 -pc_air_processor_agglom 0 -pc_air_poly_order 6 -pc_air_inverse_type arnoldi -pc_air_matrix_free_polys -pc_air_smooth_type f -pc_air_coarsest_poly_order 100 -pc_air_coarsest_inverse_type newton -pc_air_coarsest_matrix_free_polys -pc_air_auto_truncate_tol 1e-1
\end{lstlisting}
The only option we adjust in our weak scaling studies is the level at which we start testing tentative coarse grid solvers for automatic truncation of the hierarchy, as discussed in \secref{sec:truncation}. For example, for the 2 nodes shown in \fref{fig:weak_scaling_all} we start testing truncation from level 9 with \lstinline{-pc_air_auto_truncate_start_level 9}. We then increase this parameter when weak scaling based on an upper bound on the number of expected levels, given the coarsening rate in these problems. 

All of the results presented below were run on Lumi-G, a HPE Cray EX supercomputer with 2978 nodes consisting of 4 AMD MI250x GPUs and a single 64 core AMD EPYC ``Trento" CPU per node. Each of the MI250x GPUs can be considered as two GPUs linked by a fast interconnect and so we use 8 MPI ranks per node. Lumi-G has a measured Linpack performance of 380 PFLOPS and was number 9 on the June 2025 Top 500 list. We use as many DOFs as possible per rank in an attempt to keep the GPUs saturated, unless otherwise noted we use 104M DOFs per node (13M per rank) in our weak scaling studies. 
%~~~~~~~~~~~~~~
\subsection{Coarse grid solver polynomial order}
\label{sec:coarse_order}
% ~~~~~~~~~~~~~
% ~~~~~~~~~~~~~
\begin{figure}[thp]
\centering
\subfloat[][13M DOFs on 1 GPU]{\label{fig:order_serial}\includegraphics[width =0.4\textwidth]{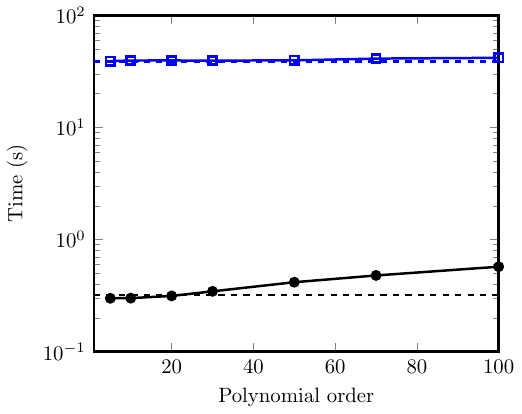}} \quad
\subfloat[][104M DOFs on 2 nodes]{\label{fig:order_2_node}\includegraphics[width =0.4\textwidth]{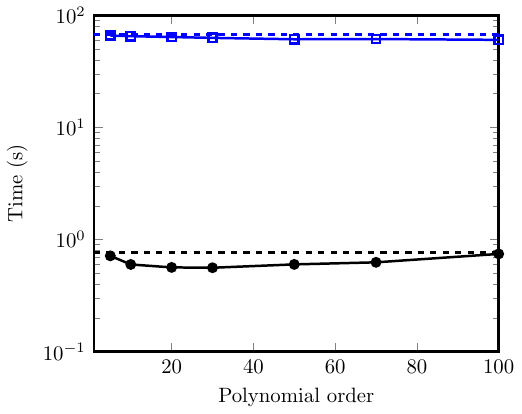}} \quad
\caption{Time vs coarse grid solver polynomial order with automatic truncation for AIRG on Lumi-G on a 2D upwind finite-difference advection problem with $\mat{v}=(\cos(\pi/4), \sin(\pi/4))$ on a structured grid. $\CIRCLE$ is the solve time, \textcolor{blue}{$\square$} is the setup. Dashed lines are the time without truncation.}
\label{fig:order}
\end{figure}
% ~~~~~~~~~~~~~
Before investigating the strong and weak scaling, we first must determine, $m_\textrm{coarse}$, the order of the polynomial we apply as the coarse grid solver, given the discussion in Sections \ref{sec:coarse grid} \& \ref{sec:truncation}. The optimal value will depend on the hierarchy generated and the balance of work to communication on a given machine, but we find even low-order default values can improve performance. \fref{fig:order} shows the effect of enabling automatic truncation on the setup and solve time given different coarse grid solver polynomial orders. 

In serial, \fref{fig:order_serial} shows that the runtime of the setup is at a minimum with no truncation and the solve is at a minimum with a low-order coarse grid solver polynomial, between $5^\textrm{th}$ and $10^\textrm{th}$ order. The solve with the $10^\textrm{th}$ order coarse grid solver polynomial is around 7\% faster than without truncation (both take 7 iterations), despite the cycle complexity increasing from 27.2 to 28.0. So although we have more work to do with truncation, we have eliminated the cost of kernel launches on many coarse grids: without truncation we have 22 levels, but with truncation we have only 13. The amount of memory used also decreases slightly, with the storage complexity at 10.22 with truncation and 10.21 without. 

In parallel on 2 nodes, \fref{fig:order_2_node} shows the benefit of truncation becomes more pronounced, which is as expected given it can help alleviate coarse grid communication bottlenecks described above. \fref{fig:order_2_node} shows that the setup time is at a minimum at $100^\textrm{th}$ order (decreasing by around 10\%), with the solve time at a minimum at $30^\textrm{th}$ order (decreasing by around 37\%). There are 25 levels, a cycle complexity of 28.5 and storage complexity of 10.5 with no truncation, and 13 levels, cycle complexity of 32.8 and storage complexity of 10.44 with $30^\textrm{th}$ order. The interplay of work to communication in the solve is particularly evident with the $100^\textrm{th}$ order polynomial, as the cycle complexity is 61.6 with 11 levels. This is $2.3\times$ the amount of work compared to the case without truncation, but results in a similar solve time. 

The reason the setup time still decreases as we increase the polynomial order is because the truncation becomes more aggressive and hence we avoid components of the setup which scale poorly on some of the coarse levels. We found in \cite{Dargaville2025} (see Figure 5) that it was the matrix-matrix products used to calculate the restrictor and the coarse grid matrix that scaled poorly in the middle of the hierarchy, due to the decreasing ratio of local to non-local work. Previously we used the combination of truncation, processor agglomeration, which decreases the number of active MPI ranks, and repartitioning of the coarse grid matrices with ParMETIS to combat this. This resulted in a large decrease in both setup and solve time, although the ParMETIS repartitioning became the largest component of our setup time. 

\begin{figure}[thp]
\centering
\subfloat[][No truncation. Total setup time: 90.2s]{\label{fig:setup_times}\includegraphics[width =0.4\textwidth]{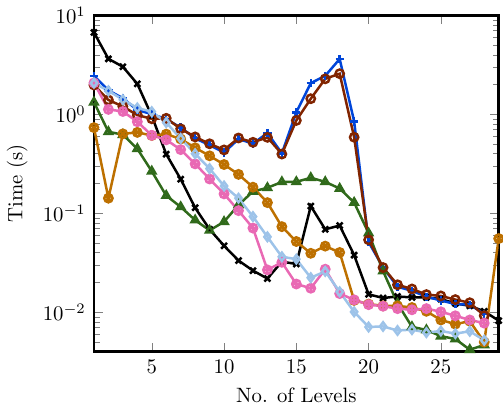}} \quad
\subfloat[][With truncation and a 100$\textrm{th}$ order polynomial as the coarse grid solver, \textcolor{gaylordpurple}{$\square$} is the truncation time. Total setup time: 68.6s]{\label{fig:setup_times_trunc}\includegraphics[width =0.4\textwidth]{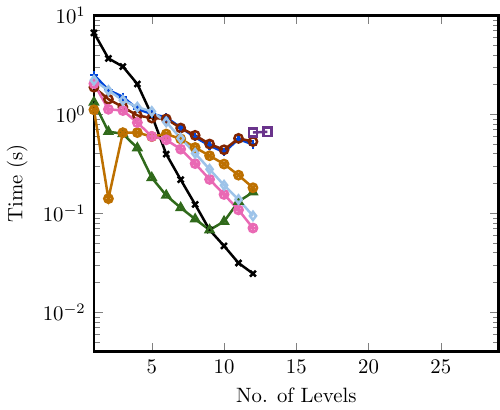}} \quad
\caption{Setup time for AIRG on Lumi-G on a 2D upwind finite-difference advection problem with $\mat{v}=(\cos(\pi/4), \sin(\pi/4))$ on a structured grid on 32 nodes with 104M DOFs/node. The \textcolor{black}{$\times$} is the CF splitting, the \textcolor{foliagegreen}{$\triangle$} is the prolongator, the \textcolor{deludedorange}{$\otimes$} is the GMRES polynomial, the \textcolor{matlabblue}{$+$} is the SpGEMM for the restrictor, the \textcolor{fireenginered}{o} is the SpGEMM for the coarse grid, the \textcolor{darksalmon}{$\oplus$} is the matrix extract, the \textcolor{skybluetext}{$\diamond$} is the dropping.}
\label{fig:setup_times_whole}
\end{figure}

We see similar growth in the time spent in the matrix-matrix products in the middle of the hierarchy, but we do not find it necessary to apply processor agglomeration in this work. The difference in approaches here is due to the combination of more DOFs per rank, which delays the level at which communication bottlenecks appear and the ability of the high-order polynomials to allow more aggressive truncation. In \cite{Dargaville2025} we only had 8.8k DOFs/rank, whereas in this work we use 13M DOFs/rank to help saturate the GPUs. For example, \fref{fig:setup_times_whole} shows the individual components of the setup time on 32 nodes. We can see in \fref{fig:setup_times} that the matrix-matrix products have a significant peak in cost starting at level 15; we find that the relative size of these peaks grows considerably as we weak scale. 

\fref{fig:setup_times_trunc} shows that with a 100$\textrm{th}$ polynomial, we can truncate so aggressively in this problem as to completely avoid these peaks. This causes the total setup time to decrease from 90.2 seconds to 68.6. The cost of this truncation is also small, as we can \textit{a priori} estimate that the truncation may first occur from level 12, and hence we only calculate the coarse grid solver and test it's efficacy on levels 12 and 13, where the truncation then occurs. The memory required to compute the polynomial coefficients in the setup is also small, as the 100 vectors required on level 13 to store the Krylov basis, for example, is equivalent to $\sim8$ vectors on the top grid. If we fix the polynomial order of the coarse grid solver, as we weak scale we will have to add levels and hence this peak will eventually reappear. Given the weak scaling results below, we found truncation with a 100$\textrm{th}$ polynomial was sufficient to avoid the worst effects up to large fractions of Lumi-G.

If we wished to scale to larger machines, we have two (complementary) options to help control this behaviour. The first is to increase the polynomial order and hence truncate more aggressively. This will create more work in the solve, but as we have seen above this may not have as large an impact as the cycle complexity suggests. The other option is to enable processor agglomeration/repartitioning as in \cite{Dargaville2025}. This can be quite expensive and there is a lack of readily available GPU repartitioning codes, but the combination of using many DOFs per rank on GPUs and truncation may limit the number of levels which require processor agglomeration/repartitioning. For the results below, we always use a 100$^\textrm{th}$ order GMRES polynomial ($m_\textrm{coarse}=100$) applied matrix-free as our coarse grid solver and we do not use processor agglomeration/repartitioning. 
% ~~~~~~~~~~~~~
% ~~~~~~~~~~~~~
%~~~~~~~~~~~~~~
\subsection{Strong scaling}
\label{sec:strong_scaling}
\begin{figure}[thp]
\centering
\subfloat[][208M DOFs from 2 nodes to 64 nodes]{\label{fig:strong_scaling}\includegraphics[width =0.4\textwidth]{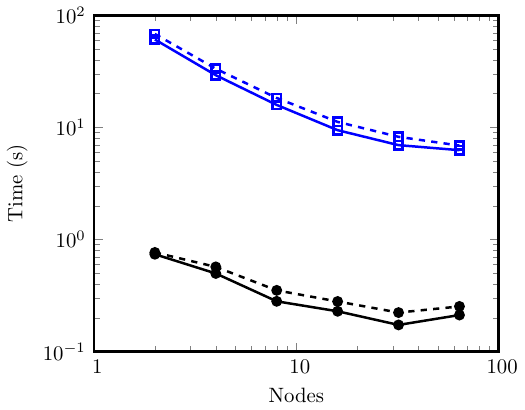}} \quad
\subfloat[][3.3B DOFs from 32 nodes to 1024 nodes]{\label{fig:strong_scaling_32_nodes}\includegraphics[width =0.4\textwidth]{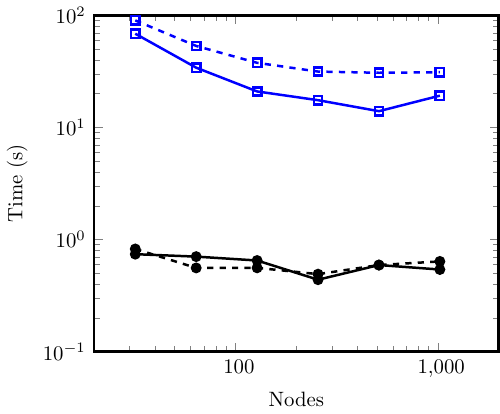}} \quad
\caption{Strong scaling results for AIRG on Lumi-G on a 2D upwind finite-difference advection problem with $\mat{v}=(\cos(\pi/4), \sin(\pi/4))$ on a structured grid. $\CIRCLE$ is the solve time, \textcolor{blue}{$\square$} is the setup. The dashed lines are without truncation.}
\label{fig:strong_scaling_all}
\end{figure}

% ~~~~~~~~~~~~~
\fref{fig:strong_scaling_all} shows the results of strong scaling studies for two different sized problems, one with 208M DOFs from 2 nodes to 64 nodes (\fref{fig:strong_scaling}) and the other with 3.3B DOFs from 32 nodes to 1024 nodes (\fref{fig:strong_scaling_32_nodes}). We start both scaling studies with 104M DOFs/node, which is close to the maximum we can fit per node in this problem.

We see that the setup time for the smaller problem in \fref{fig:strong_scaling} continues to decrease out to 64 nodes, giving a strong scaling efficiency of 30\%. The solve reaches a minimum around 32 nodes with a strong scaling efficiency of 27\% and a throughput of 2.6\xten{8} DOFs/s/node/V-cycle, down from 9.8\xten{8} DOFs/s/node/V-cycle with 2 nodes. The iteration count has remained the same (at seven iterations) but the performance has decreased, which is to be expected as GPUs perform best when fully saturated. The 32 node case for example has only 812k DOFs/rank, compared with 13M DOFs/rank with 2 nodes. We can see that the performance suffers in both the setup and solve if we do not truncate as we strong scale. 

For the larger problem in \fref{fig:strong_scaling_32_nodes}, our strong scaling efficiency in the solve decreases when compared with the smaller problem, with the solve time reaching a minimum at 256 nodes with an efficiency of 21\%. We also see that the truncation is not always helpful in the solve for the larger problem. The setup time however strong scales similarly to the smaller problem, reaching a minimum at 512 nodes with efficiency of 31\%. Interestingly, we see that the truncation has improved the performance of the setup considerably. We saw this in \fref{fig:setup_times_whole} (which is the 32 node case in \fref{fig:strong_scaling_32_nodes}), where the setup decreased from 90.2s to 68.6s with truncation (a factor of 1.3$\times$). This factor increases as we strong scale, for example at 512 nodes the setup is 2.2$\times$ faster with truncation than without. This matches the discussion in \secref{sec:coarse_order} and in \cite{Dargaville2025}; we find that the peaks in setup time caused by the matrix-matrix products increase as the ratio of local to non-local work decreases, which occurs as we strong scale. 

\fref{fig:strong_scaling_32_nodes} shows that the solve time initially decreases with truncation, before increasing slightly until around 512 nodes before decreasing again after 512 nodes. We can make the general statement that as we strong scale, we eventually reach the point where there is not enough work on the lower grids to hide the communication and this is where truncation becomes important in the solve in the larger problem. The same is true when we weak scale, as there will be more levels in the hierarchy. 

Predicting when this balance pays off in the solve is difficult for any one problem e.g., we see consistent reduction in solve time in the smaller problem but not the larger. We know that with truncation we have fewer levels in the hierarchy and hence the coarse grid solver does more work. To see a decrease in solve time with truncation this increased work must require less time than performing the work on the truncated coarse levels, which are communication dominated. Our observation is that the solve is more resistant to fewer DOFs/rank than the setup. The only operation occurring in the solve is SpMVs (with asynchronous communication) on all levels and we find these strong scale on lower levels better than many components of the setup. 

In \cite{Dargaville2025} we operated in a regime where communication bottlenecks dominated on the lower levels, even in the solve, as we were using only 8.8k DOFs/rank on CPUs. Hence we saw consistent improvement in solve time with truncation (and processor agglomeration/repartitioning). Building a performance model would be helpful in trying to predict this behaviour, but can be difficult in algebraic multigrid methods given the interplay of factors such as increasingly non-local coarse grid matrices. 

These strong scaling studies show that the highest throughput is achieved with many DOFs/rank (13M), but decreases in wall times can be achieved in the solve out to 800k--1.6M per rank. The setup times continue to decrease in both the small and large problems while strong scaling and we see that truncation can play a role in decreasing the setup and solve times at scale. For the weak scaling studies in the next section, we use 13M DOFs/rank (104M DOFs/node) unless otherwise noted and always apply truncation.
% ~~~~~~~~~~~~~
% ~~~~~~~~~~~~~
%~~~~~~~~~~~~~~
\subsection{Weak scaling}
\label{sec:weak_scaling_square}
% ~~~~~~~~~~~~~
% ~~~~~~~~~~~~~
\begin{figure}[thp]
\centering
\subfloat[][$\CIRCLE$ is the solve time, \textcolor{blue}{$\square$} is the setup. AIRG are the solid lines with nAIR the dotted.]{\label{fig:weak_scaling}\includegraphics[width =0.4\textwidth]{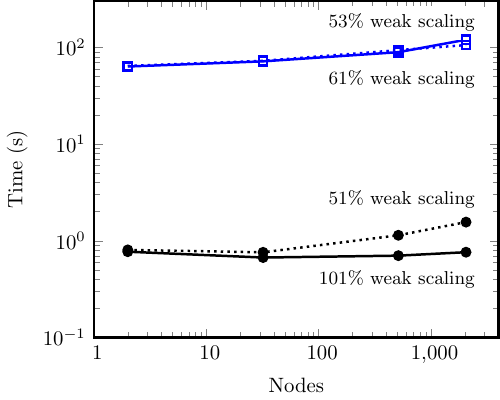}} \quad
\subfloat[][Components of the AIRG setup. The \textcolor{black}{$\times$} is the CF splitting, the \textcolor{foliagegreen}{$\triangle$} is the prolongator, the \textcolor{deludedorange}{$\otimes$} is the GMRES polynomial, the \textcolor{matlabblue}{$+$} is the SpGEMM for the restrictor, the \textcolor{fireenginered}{o} is the SpGEMM for the coarse grid, the \textcolor{darksalmon}{$\oplus$} is the matrix extract, the \textcolor{skybluetext}{$\diamond$} is the dropping and \textcolor{gaylordpurple}{$\square$} is the truncation.]{\label{fig:setup_scaling}\includegraphics[width =0.4\textwidth]{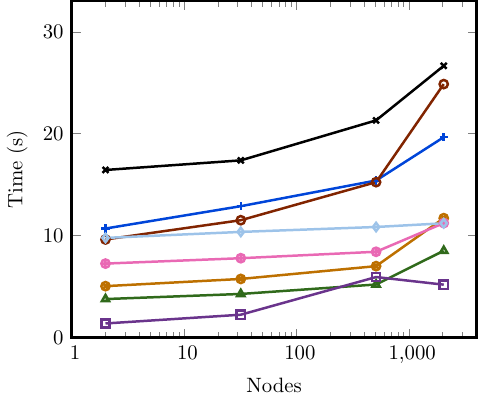}} \quad
\caption{Weak scaling results for AIRG on Lumi-G on a 2D upwind finite-difference advection problem with $\mat{v}=(\cos(\pi/4), \sin(\pi/4))$ on a structured quad mesh, with 104M DOFs/node from 2 nodes (208M DOFs) to 2048 nodes (213B DOFs).}
\label{fig:weak_scaling_all}
\end{figure}
\begin{figure}[thp]
\centering
\subfloat[][Number of non-zeros in the coarse matrix on each level relative to the top grid matrix. $\CIRCLE$ is on 2 nodes, \textcolor{blue}{$\square$} is on 32 nodes.]{\label{fig:rel_coarse_grid_nnzs}\includegraphics[width =0.4\textwidth]{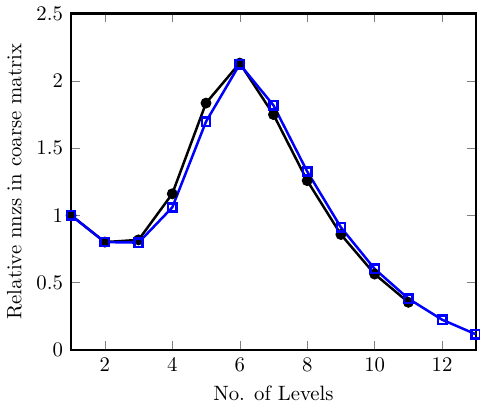}} \quad
\subfloat[][Cycle complexity during weak scaling with truncation]{\label{fig:weak_scaling_discontinuity}\includegraphics[width =0.4\textwidth]{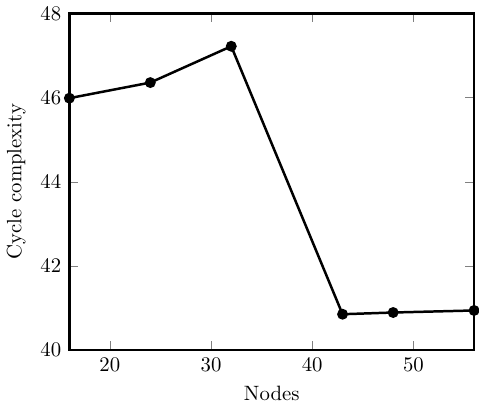}} \quad
\caption{Effect of truncation on the amount of work during weak scaling.}
\label{fig:tempy}
\end{figure}

% ~~~~~~~~~~~~~
\fref{fig:weak_scaling_all} shows weak scaling results, from 2 nodes with 208M DOFs (14,422$^2$ mesh) to 2048 nodes with 213B DOFs (461,511$^2$ mesh). In this example we increase the number of DDC iterations from 2 to 3 compared with the previous section (\lstinline{-pc_air_ddc_its 3}) in order to improve the convergence slightly. We see that AIRG weak scales very well in this problem, with 101\% weak scaling efficiency in the solve and 53\% in the setup. The iteration count in the solve grows from 7 to 8; despite this we find greater than ideal weak scaling. This is due to the use of truncation. \fref{fig:rel_coarse_grid_nnzs} shows the number of non-zeros in the coarse matrices on each level relative to the top grid matrix, for both 2 nodes and 32 nodes. We can see that as we weak scale, the relative density of our operators does not change. As noted above, we use a 100$^\textrm{th}$ order GMRES polynomial as our coarse grid solver and apply automatic truncation. For the 2 node problem, we can see in \fref{fig:rel_coarse_grid_nnzs} this results in 11 levels in the hierarchy, with 13 levels for the 32 node case. The relative number of non-zeros in the matrix on level 11 is around 0.38, whereas at level 13 it is 0.116. The coarse grid solver doing 100 SpMVs on level 11 therefore costs around 38 top grid SpMVs, compared with 11.6 on level 13. 
% ~~~~~~~~~~~~~
% ~~~~~~~~~~~~~
As we weak scale with an iterative coarse grid solver and a truncated hierarchy, the relative amount of work our coarse grid solver performs therefore decreases. This is evident in the cycle complexity, which for the weak scaling shown is 74, 43, 35.8 and 33.5, plateauing with increasing resolution. We also see the same with heavy truncation in Table 8 from \cite{Dargaville2025}. 

We should note that this doesn't mean we are artificially increasing the time taken to solve at smaller scales. For example, both \fref{fig:order_2_node} and \fref{fig:strong_scaling} show that with a 100$^\textrm{th}$ order GMRES polynomial as a coarse grid solver and automatic truncation we have the same solve time for 2 nodes, despite the cycle complexity increasing. \fref{fig:strong_scaling_32_nodes} shows the same is true with 32 nodes, with the solve time with truncation actually decreasing slightly. 

As we weak scale and add more levels, the relative size of the work performed on the coarse grid solve decreases and we see the cycle complexity with truncation approach that without. This is not an entirely smooth transition however, as the amount of work increases slightly before we add a new level, at which point there is a large discontinuous decrease, before beginning to increase again. We can see that in \fref{fig:weak_scaling_discontinuity}, which focuses on the discontinuity that occurs between the 32 and 43 node case, where the number of levels increases from 13 to 14. Hence depending on where each point in our weak scaling study lands, we can see small variations in the solve time. This is what causes the slight improvement in \fref{fig:weak_scaling}; we do not claim that our methods are weak scaling better than ideal in the limit. This effect would also be evident in other multigrids that use an iterative coarse grid solver and truncation. 

The storage complexity in this problem stays very close to constant, at 10.5, 10.9, 10.98 and 10.96 as we weak scale. This means our memory use, although high when compared to an elliptic multigrid, does not grow substantially. We note that this is almost twice that shown in \cite{Dargaville2025}, where the storage complexity is between 5--6. This is the impact of using structured meshes, they require slower coarsenings with more levels than unstructured meshes in advection problems, as noted in \secref{sec:aim}. Furthermore, the memory use in this problem could be almost twice as large if not for two choices we described above, namely performing only F-point smoothing (described in \secref{sec:reduction}) and using the GMRES polynomials matrix-free when smoothing (described in \secref{sec:airg}). For example, on 2 nodes, as noted, we store the equivalent of 10.5 copies of the top grid matrix. If we add in matrix-free C-point smoothing, this becomes 17.2, due to the need to store $\mat{A}_\textrm{cc}$ and $\mat{A}_\textrm{cf}$. If we then also use an assembled approximation to our polynomials to perform F and C point smoothing, i.e., $\bar{\mat{A}}_\textrm{ff}^{-1} = \hat{\mat{A}}_\textrm{ff}^{-1}$, this increases to 20.2. This helps show the power of the GMRES polynomials, as they allow matrix-free application while also being strong enough to not require the use of C-point smoothing. 

\fref{fig:setup_scaling} shows the weak scaling of each component of our setup for AIRG. We can see that all aspects of the setup take more time as we scale, with the biggest components the matrix-matrix products and the CF splitting. Typically the matrix-matrix products dominate multigrid setup, but our asymmetric linear system requires us to perform a CF splitting with $\mat{S}+\mat{S}^\textrm{T}$ and this is more expensive than typical in elliptic multigrids. The weak scaling we see in our setup is much better than in \cite{Dargaville2025}, where we saw between 8--19\% efficiency. Again this is due to using more DOFs/rank in this work.
%~~~~~~~~~~~~~~
\subsubsection{Comparison with nAIR}
\label{sec:nair}
% ~~~~~~~~~~~~~
It is difficult to perform comparisons with other reduction multigrids on GPUs, as we believe that the only GPU port of a reduction multigrid currently available outside of \textsc{PFLARE} is the nAIR \cite{Manteuffel2019} implementation in \textit{hypre}. We do not believe this would be a fair comparison for several reasons, namely:
\begin{enumerate}
\item The only GPU CF splitting available in \textit{hypre} is PMIS and we showed in \cite{Dargaville2025} that this gave very poor convergence in reduction multigrids on asymmetric problems.
\item \textit{hypre} cannot use the Neumann polynomials as F-point smoothers; the literature for nAIR uses Jacobi smoothing instead, which we showed in \cite{Dargaville2025} requires more work in the solve.
\item \textit{hypre} assembles matrix representations of the Neumann polynomials to create an approximate ideal restrictor with a matrix-matrix product, but does not use fixed sparsity like in \eref{eq:fixed_sparsity}, making the use of polynomials of order $>1$ very expensive. 
\item There are limited choices for iterative coarse-grid solvers and truncation in \textit{hypre} and hence we would lack the ability to remove coarse grid bottlenecks.
\end{enumerate}

Given this, we built a GPU implementation of nAIR in PFLARE that can apply the Neumann polynomials matrix-free for use with F-point smoothing and build assembled polynomial approximations with fixed sparsity to compute an approximate ideal restrictor. We therefore use the same PMISR DDC splitting, polynomial order and coarse grid solver options with nAIR as AIRG. The only change required to the command-line options is modifying \lstinline{-pc_air_inverse_type arnoldi} to \lstinline{-pc_air_inverse_type neumann}. This helps us evaluate the impact of our GMRES polynomials separate from other implementation details for a fair comparison. 

We can see in \fref{fig:weak_scaling} that the setup of nAIR is slightly cheaper than AIRG and scales better, at 61\% compared with 53\% for AIRG. The solve however, only weak scales with 51\% efficiency. This is due to poorer convergence, with the iteration count increasing from 7 with 2 nodes to 15 with 2048 nodes. The time per iteration is very similar to AIRG (given the same polynomial orders are used), but nAIR is not scalable in this problem due to the growing iteration count. Given this we do not investigate the use of nAIR any further.
%~~~~~~~~~~~~~~
\subsubsection{Dependence on direction}
\label{sec:dependence}
% ~~~~~~~~~~~~~
We noted in \secref{sec:aim} that the convergence of reduction multigrids on structured meshes can depend on the velocity direction in advection problems. This is because the size of the largest diagonally dominant $\mat{A}_\textrm{ff}$ submatrix that can be extracted depends on direction, given the size of the matrix entries changes. This is one disadvantage to using a reduction multigrid on a structured grid when compared with a sweep. There is no dependence on direction in the convergence of a sweep, given each row/element block is exactly inverted. 

In the previous section, we specified $\mat{v}=(\cos(\pi/4), \sin(\pi/4))$ and showed that not only could we get good convergence with AIRG, but that the resulting solve time is scalable. We note however that this is one of the easiest problems to solve on a structured mesh (see Figure 4 in \cite{Manteuffel2019}). The grid-aligned directions, $\mat{v}=(1, 0)$ and $\mat{v}=(0, 1)$ are also very easy to solve, given this would reduce to a set of uncoupled 1D advection equations. \secref{sec:reduction} noted that 1D advection problems fit into a tridiagonal stencil and hence a reduction multigrid is equivalent to a cyclic reduction method and can become an exact solver.

\begin{figure}[thp]
\centering
\subfloat[][$\CIRCLE$ is the solve time, \textcolor{blue}{$\square$} is the setup.]{\label{fig:weak_scaling_direction_sqrt}\includegraphics[width =0.4\textwidth]{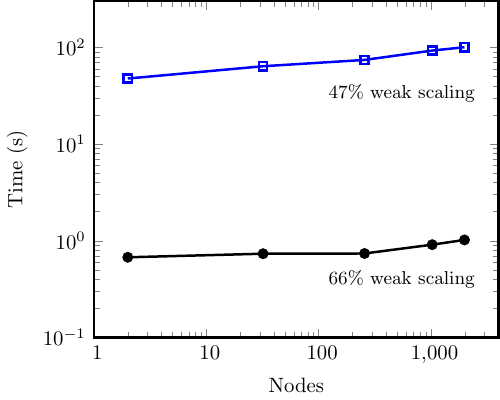}} \quad
\subfloat[][Components of the setup. The \textcolor{black}{$\times$} is the CF splitting, the \textcolor{foliagegreen}{$\triangle$} is the prolongator, the \textcolor{deludedorange}{$\otimes$} is the GMRES polynomial, the \textcolor{matlabblue}{$+$} is the SpGEMM for the restrictor, the \textcolor{fireenginered}{o} is the SpGEMM for the coarse grid, the \textcolor{darksalmon}{$\oplus$} is the matrix extract, the \textcolor{skybluetext}{$\diamond$} is the dropping and \textcolor{gaylordpurple}{$\square$} is the truncation.]{\label{fig:setup_scaling_direction_sqrt}\includegraphics[width =0.4\textwidth]{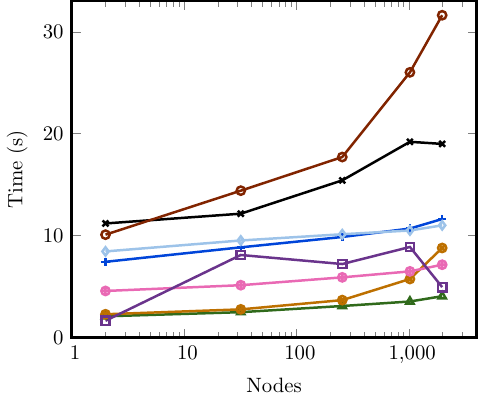}} \quad
\caption{Weak scaling results for AIRG on Lumi-G on a 2D upwind finite-difference advection problem with $\mat{v}=(\sqrt{2/3}, \sqrt{1/3})$ on a structured quad mesh, with 64M DOFs/node from 2 nodes (128M DOFs) to 1984 nodes (127B DOFs).}
\label{fig:weak_scaling_direction_sqrt_all}
\end{figure}

One of the most challenging directions to solve is $\mat{v}=(\sqrt{2/3}, \sqrt{1/3})$ and hence \fref{fig:weak_scaling_direction_sqrt_all} shows the result from a weak scaling study with this velocity. The fact that some directions are harder to solve than others does not mean we can't achieve good convergence and/or scalability across different directions. We find for this more difficult direction we have to modify solver parameters to achieve good performance. We could have used the same modified parameters in the easier $\pi/4$ problem and still seen scalable convergence, but this would be more expensive than tailoring the solver parameters for the difficulty of each problem.  

In particular, we have to slow the coarsening compared to the $\pi/4$ case in order to retain good diagonal dominance in $\mat{A}_\textrm{ff}$. The only change required to the command-line options shown above in \secref{sec:Results} is modifying \lstinline{-pc_air_strong_threshold 0.99} to \lstinline{-pc_air_strong_threshold 0.4}. Given this we also have to decrease the number of DOFs/rank to 8M (64M DOFs/node), as the slower coarsening requires more levels and hence more memory. The grid complexity increases to around 3.2, compared with 2.8 for the $\pi/4$ direction. 
% ~~~~~~~~~~~~~
% ~~~~~~~~~~~~~
We can see in \fref{fig:weak_scaling_direction_sqrt_all} that both the solve and setup weak scale with lower efficiencies than in \fref{fig:weak_scaling_all}. We still retain 66\% weak scaling efficiency, and this growth is due to both an increasing iteration count, from 6 to 8 when moving from 2 nodes to 1984 nodes and the greater number of levels. Given the slower coarsening, the number of levels increases from 16 to 35 and this also drives the poorer weak scaling seen in the setup. We can see in \fref{fig:setup_scaling_direction_sqrt} that the time required to compute the coarse-grid matrix-matrix products dominates, in contrast to \fref{fig:setup_scaling} which required only 17 levels with 2048 nodes. The cycle and storage complexities are also higher, with 106, 61, 47.2, 39.5 and 37.4 for the cycle and 12.2, 12.8, 13.0, 13.1 and 13.1 for the storage, respectively. The throughput is roughly half of that seen with the $\pi/4$ case, at approximately 5\xten{8} DOFs/s/node/V-cycle. 

Although the performance clearly depends on direction, we can still obtain excellent convergence and good weak scaling behaviour in structured grid advection problems with different directions. As noted in \secref{sec:reduction}, the fact that reduction multigrids can become exact solvers as the approximation to $\mat{A}_\textrm{ff}^{-1}$ is improved is key to this behaviour. In order to minimise the number of parameters we changed in this problem, we only decreased the coarsening rate, but there are a range of parameters we can modify to improve the convergence in asymmetric problems. In rough order of importance, we can:
\begin{enumerate}
\item Slow the coarsening rate, by decreasing the strong threshold in the PMISR (\lstinline{-pc_air_strong_threshold}) and/or increasing the DDC fraction or number of iterations (\lstinline{-pc_air_ddc_fraction -pc_air_ddc_its})
\item Decrease the drop tolerances on $\mat{A}$ and/or $\mat{R}$ on each level (\lstinline{-pc_air_a_drop -pc_air_r_drop})
\item Increase the GMRES polynomial order (\lstinline{-pc_air_poly_order})
\item Add C-point smoothing (\lstinline{-pc_air_smooth_type fc}) or multiple iterations of either \\ e.g., FCF smoothing (\lstinline{-pc_air_smooth_type fcf})
\item Increase the nnzs retained in $\hat{\mat{A}}_\textrm{ff}^{-1}$ when computing the approximate ideal operators  \\(\lstinline{-pc_air_inverse_sparsity_order})
\item Use an approximate ideal prolongator (\lstinline{-pc_air_one_point_classical_prolong 0})
\item Improve the approximate ideal operators with a Richardson iteration \\(\lstinline{-pc_air_improve_z_its 1 -pc_air_improve_w_its 1})
\end{enumerate}
The impact of each of these parameters is problem dependent and the types of asymmetric problems with which we can obtain scalable convergence is an open question, but we have found good performance with reduction multigrids in many difficult, asymmetric problems.
% ~~~~~~~~~~~~~
%~~~~~~~~~~~~~~
\subsection{Weak scaling - rectangular}
\label{sec:weak_scaling_rectangular}
% ~~~~~~~~~~~~~
Given we have now shown good convergence and scaling in square problems, we decided to test the performance of our method in a rectangular problem. Rectangular problems can be well suited to sweeps using KBA decompositions (\fref{fig:pencil}) in parallel. In 2D, for good performance in a sweep we would require enough resolution in one dimension (say the $x$ direction) to both have enough pencils to distribute one to all ranks and ensure that each pencil has sufficient local work (i.e., wide pencils). If those criteria are satisfied, we then require as much resolution as possible in the $y$ dimension. This ensures the relative number of start-up stages is small (i.e., the pipefill) and hence we would expect to see excellent parallel utilisation as almost every rank would be working throughout the entire computation. 

If we were then performing a weak scaling study and increased the resolution in only the $x$ dimension we would expect ideal weak scaling efficiency (assuming there is enough $y$ resolution in the original problem to keep the relative impact of the pipefill low), as we could distribute the extra pencils to each of the extra ranks. If instead we increased the resolution only in the $y$ direction, this would scale poorly as we would be forced to start decomposing in $y$, effectively serializing the sweep in such a rectangular problem. We decided to test the performance of our reduction multigrid in this limit, as it would be pathological for a sweep algorithm.

We start with 2 nodes (16 ranks) solving on a rectangular domain, $L_x=1, L_y=832$, with rectangular resolution, namely a $500 \times 416,000$ mesh. With a KBA decomposition, we could create $16 \times 1$ pencils of size approximately $31 \times 416,000$ and on 2 nodes we would likely see good parallel utilisation from a sweep. As we weak scale however, we increase both the domain size and resolution only in $y$, with the $x$ dimension domain size and resolution remaining fixed. The reason for modifying both the domain and resolution is to ensure we retain the same FD stencil. 

When weak scaling out to 2800 nodes (22,400 ranks) this results in a rectangular domain with $L_x=1, L_y=\text{1,164,800}$ and a $500 \times \text{582,400,000}$ mesh. There is no way to create enough pencils to fill all the ranks with a KBA decomposition without decomposing in the $y$ direction. On GPUs a KBA sweep would also be forced to balance shared memory and distributed parallelism. For example, on 2800 nodes, if we decomposed in both $x$ and $y$ using a KBA decomposition with $500 \times 44$ pencils in an attempt to maximize the distributed parallelism, each pencil would have a size of around $1 \times 13\textrm{M}$, leaving no parallelism available for each GPU in the local sweep. Our goal in this example is to show we can achieve good performance in a problem where a sweep could not. As in the previous examples we use a volumetric decomposition with our reduction multigrid.
% ~~~~~~~~~~~~~
\begin{figure}[thp]
\centering
\subfloat[][$\CIRCLE$ is the solve time, \textcolor{blue}{$\square$} is the setup.]{\label{fig:weak_scaling_rectangular}\includegraphics[width =0.4\textwidth]{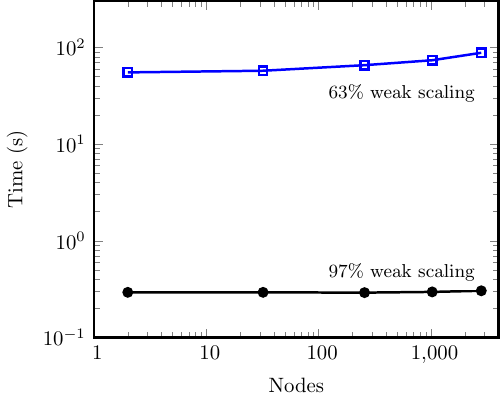}} \quad
\subfloat[][Components of the setup. The \textcolor{black}{$\times$} is the CF splitting, the \textcolor{foliagegreen}{$\triangle$} is the prolongator, the \textcolor{deludedorange}{$\otimes$} is the GMRES polynomial, the \textcolor{matlabblue}{$+$} is the SpGEMM for the restrictor, the \textcolor{fireenginered}{o} is the SpGEMM for the coarse grid, the \textcolor{darksalmon}{$\oplus$} is the matrix extract, the \textcolor{skybluetext}{$\diamond$} is the dropping and \textcolor{gaylordpurple}{$\square$} is the truncation.]{\label{fig:setup_scaling_rectangular}\includegraphics[width =0.4\textwidth]{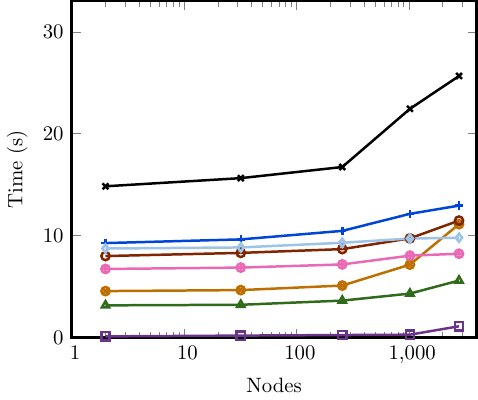}} \quad
\caption{Weak scaling results for AIRG on Lumi-G on a 2D upwind finite-difference advection problem with $\mat{v}=(\cos(\pi/4), \sin(\pi/4))$ on a structured quad mesh, with 104M DOFs/node from 2 nodes (208M DOFs) to 2800 nodes (291B DOFs). The domain size and resolution are increased only in the $y$ dimension.}
\label{fig:weak_scaling_rectangular_all}
\end{figure}
% ~~~~~~~~~~~~~

We use $\mat{v}=(\cos(\pi/4), \sin(\pi/4))$ in this problem and note that because information does not propagate the length of the longest dimension, this problem is easier to solve for iterative methods than a square problem. Given the structured mesh, changing the direction to $\mat{v}=(0, 1)$ however would not make the problem more challenging for a reduction multigrid, given the grid aligned direction could result in an exact solver (see the discussion in \secref{sec:dependence}). We exploit the easier problem and modify one of the solver parameters to improve the performance, similar to how we tailored parameters to specific problems in \secref{sec:dependence}. In particular, we decrease the coarse grid solver polynomial order to 10. The only change required to the command-line options shown above in \secref{sec:Results} is modifying \lstinline{-pc_air_coarsest_poly_order 100} to \lstinline{-pc_air_coarsest_poly_order 10}. We found a 100$^\textrm{th}$ GMRES polynomial was overkill as a coarse grid solver in this problem, as a 10$^\textrm{th}$ order polynomial with automatic truncation results in a constant 10 levels in our hierarchy as we weak scale. 

\fref{fig:weak_scaling_rectangular} shows that AIRG achieves 97\% weak scaling efficiency in the solve and 63\% in the setup for this rectangular problem. The iteration count is constant at 6 iterations and we find a constant cycle and storage complexity at approximately 28.2 and 9.5, respectively. This has the effect of improving the throughput when compared to the square problems, we achieve 2 \xten{9} DOFs/s/node/V-cycle. 
%~~~~~~~~~~~~~~~~~~~~~~~~~~~~~~~~~~~~
%~~~~~~~~~~~~~~~~~~~~~~~~~~~~~~~~~~~~
\section{Conclusions}
We believe we are the first to achieve good parallel weak scaling efficiency in the solve for single direction advection equations. We used the combination of the reduction multigrid, AIRG, and the CF splitting, PMISR DDC, to solve 2D upwind FD discretisations of advection equations on both square and rectangular domains with structured meshes. Our runs were performed on Lumi-G, a pre-exascale GPU machine, and we tailored our methods for use on GPUs. We tested our weak scaling from 2 nodes to between 1984 (67\% of Lumi-G) and 2800 nodes (94\% of Lumi-G). 

We showed results for two different velocity directions and though we found the performance was dependent on the direction with structured meshes, we were are able to obtain good weak scaling in both cases, with 66\% efficiency in the solve for $\mat{v}=(\sqrt{2/3}, \sqrt{1/3})$ up to 1984 nodes and 101\% in the $\pi/4$ case up to 2048 nodes with AIRG. We found the setup times weak scaled with 47\% and 53\% efficiency, respectively. In the $\pi/4$ case, we also tested an implementation of a different reduction multigrid, nAIR, and found much poorer weak scaling efficiency in the solve at 51\%. In a rectangular problem, we used different mesh resolution in each dimension in a way that would challenge a sweep algorithm. We found 97\% weak scaling efficiency in the solve and 63\% in the setup up to 2800 nodes with AIRG. Across all the problems tested, we achieved a throughput of between 5 \xten{8} -- 2 \xten{9} DOFs/s/node/V-cycle in the solve.

To facilitate their use the methods in this paper are available in the open-source PFLARE library. \textsc{PFLARE} provides new PC types for PETSc and can be easily used with existing PETSc code, as it is available directly through the PETSc configure (\texttt{--download-pflare}). \textsc{PFLARE} runs on CPUs and GPUs (AMD, NVIDIA or Intel) with interfaces in C/Fortran/Python. There are several advantages to using these methods when compared to other ways of solving advection-type problems, including:
\begin{enumerate}
\item Our use of a volumetric decomposition in all the examples shown, which matches the decomposition used in other multi-physics. We do not require KBA-style pencil decompositions and pipelining of multiple directions for good parallel performance.
\item We do not depend on sweeps/Gauss-Seidel methods, or any specific ordering of unknowns.
\item We can solve problems with lower-triangular structure, like the upwind FD discretisations used in this work, or without, like the stabilised CG FEM discretisation in \cite{Dargaville2025}. 
\item Convergence, memory use and parallelism often improve on unstructured meshes (due to fewer levels).
\item Most of the operations required in the solve and setup are common linear algebra kernels, such as sparse matrix-vector and matrix-matrix products which have optimised implementations on different hardware.
\item We have two mechanisms for improving parallel performance on large machines, truncation of the hierarchy (which we used in this paper) and processor agglomeration/repartitioning (which we did not). 
\end{enumerate}
Some disadvantages include:
\begin{enumerate}
\item The setup can be expensive and weak scales less well than the solve, as is common with algebraic multigrid methods. We found our setup costs between 100-290 solves, depending on the problem, though this cost can often be amortised by performing many solves (e.g., time-dependent problems).     
\item The memory use can be high. With upwind FD discretisations on structured meshes we required between 10-13 copies of the top-grid matrix. 
\end{enumerate} 
Overall we believe reduction multigrids offer new capabilities for solving asymmetric linear systems at scale and our future work will investigate using these methods in a variety of applications. 
%~~~~~~~~~~~~~~~~~~~~~~~~~~~~~~~~~~~~
\section*{Acknowledgements}
We acknowledge the EuroHPC Joint Undertaking for awarding this project access to the EuroHPC supercomputer LUMI, hosted by CSC (Finland) and the LUMI consortium through a EuroHPC Benchmark Access call. We acknowledge the support of the EPSRC through the funding of the EPSRC grant EP/T000414/1. We would also like to thank the PETSc team (particularly Junchao Zhang) for their assistance. 
%~~~~~~~~~~~~~~~~~~~~~~~~~~~~~~~~~~~~
%~~~~~~~~~~~~~~~~~~~~~~~~~~~~~~~~~~~~
%~~~~~~~~~~~~~~~~~~~~~~~~~~~~~~~~~~~~



%% References

\bibliographystyle{model1-num-names}
\bibliography{bib_library}

@Article{Adams2013,
  Title                    = {Provably optimal parallel transport sweeps on regular grids},
  Author                   = {Adams, M. P. and Adams, M. L. and Hawkins, W. D. and Smith, T. and Rauchwerger, L. and Amato, N. M. and Bailey, T. S. and Falgout, R. D.},
  Journal                  = {Proc. International Conference on Mathematics and Computational Methods Applied to Nuclear Science \& Engineering},
  Year                     = {2013},
  Language                 = {English},
  Url                      = {http://inis.iaea.org/Search/search.aspx?orig_q=RN:45033853},
  Urldate                  = {2015-01-15}
}

@InProceedings{Baker2012a,
  author     = {Baker, C. and Davidson, G. and Evans, T.M. and Hamilton, S. and Jarrell, J. and Joubert, W.},
  title      = {High performance radiation transport simulations: {Preparing} for {TITAN}},
  booktitle  = {{SC} '12: {Proceedings} of the {International} {Conference} on {High} {Performance} {Computing}, {Networking}, {Storage} and {Analysis}},
  year       = {2012},
  pages      = {1--10},
  month      = nov,
  doi        = {10.1109/SC.2012.64},
  shorttitle = {High performance radiation transport simulations},
}

@Article{Baker1998,
  Title                    = {An {Sn} {Algorithm} for the {Massively} {Parallel} {CM}-200 {Computer}},
  Author                   = {Baker, Randal and Koch, Kenneth},
  Journal                  = {Nuclear Science and Engineering},
  Year                     = {1998},

  Month                    = mar,
  Pages                    = {312--320},

  Doi                      = {10.13182/NSE98-1},
  ISSN                     = {00295639},
  Url                      = {http://epubs.ans.org/?a=1958},
  Urldate                  = {2016-03-10}
}

@Article{Colomer2013,
  Title                    = {Parallel algorithms for transport sweeps on unstructured meshes},
  Author                   = {Colomer, G. and Borrell, R. and Trias, F. X. and Rodríguez, I.},
  Journal                  = {Journal of Computational Physics},
  Year                     = {2013},

  Month                    = jan,
  Number                   = {1},
  Pages                    = {118--135},
  Volume                   = {232},
  ISSN                     = {0021-9991},
  Urldate                  = {2015-01-15}
}

@Article{karypis_metis-unstructured_1995,
  Title                    = {Metis-unstructured graph partitioning and sparse matrix ordering system, version 2.0},
  Author                   = {Karypis, George and Kumar, Vipin},
  Year                     = {1995},
  Urldate                  = {2014-09-02}
}

@Article{Nachtigal1992,
  Title                    = {A hybrid {GMRES} algorithm for nonsymmetric linear systems},
  Author                   = {Nachtigal, Noël M. and Reichel, Lothar and Trefethen, Lloyd N.},
  Journal                  = {SIAM Journal on Matrix Analysis and Applications},
  Year                     = {1992},
  Number                   = {3},
  Pages                    = {796--825},
  Volume                   = {13},
  Url                      = {http://epubs.siam.org/doi/abs/10.1137/0613050},
  Urldate                  = {2016-09-07}
}

@InProceedings{Nowak1999,
  Title                    = {Radiation transport calculations on unstructured grids using a spatially decomposed and threaded algorithm},
  Author                   = {Nowak, P. and Nemanic, M. K.},
  Booktitle                = {Proc. of {Int}. {Conf}. {Mathematics} and {Computation}, {Reactor} {Physics} and {Environmental} {Analysis} in {Nuclear} {Application}, {Madrid}, {Spain}},
  Year                     = {1999},
  Pages                    = {379--390},
  Volume                   = {1},
  Url                      = {https://e-reports-ext.llnl.gov/pdf/236483.pdf},
  Urldate                  = {2015-01-15}
}

@Article{Pautz2002,
  Title                    = {An algorithm for parallel {Sn} sweeps on unstructured meshes},
  Author                   = {Pautz, Shawn D.},
  Journal                  = {Nuclear science and engineering},
  Year                     = {2002},
  Number                   = {2},
  Pages                    = {111--136},
  Volume                   = {140},
  ISSN                     = {0029-5639},
  Language                 = {eng},
  Url                      = {http://cat.inist.fr/?aModele=afficheN&cpsidt=13479645},
  Urldate                  = {2015-01-15}
}

@Article{Plimpton2005,
  Title                    = {Parallel {Sn} {Sweeps} on {Unstructured} {Grids}: {Algorithms} for {Prioritization}, {Grid} {Partitioning}, and {Cycle} {Detection}},
  Author                   = {Plimpton, Steven J. and Hendrickson, Bruce and Burns, Shawn P. and Iii, William Mclendon},
  Journal                  = {Nuclear Science and Engineering},
  Year                     = {2005},
  Pages                    = {267--283},
  Volume                   = {150},
  Shorttitle               = {Parallel {Sn} {Sweeps} on {Unstructured} {Grids}}
}

@Article{saad_gmres:_1986,
  author     = {Saad, Youcef and Schultz, Martin H.},
  title      = {{GMRES}: A generalized minimal residual algorithm for solving nonsymmetric linear systems},
  journal    = {{SIAM} Journal on scientific and statistical computing},
  year       = {1986},
  volume     = {7},
  number     = {3},
  pages      = {856--869},
  shorttitle = {{GMRES}},
  urldate    = {2014-09-02},
}

@Article{Sterck2006,
  Title                    = {Reducing complexity in parallel algebraic multigrid preconditioners},
  Author                   = {Sterck, Hans De and Yang, Ulrike Meier and {Jeffrey} and Heys, J.},
  Journal                  = {SIAM J. Matrix Anal. Appl},
  Year                     = {2006},
  Pages                    = {1019--1039},
  Volume                   = {27}
}

@Article{Manteuffel2019,
  author   = {Manteuffel, Thomas A. and Münzenmaier, Steffen and Ruge, John and Southworth, Ben},
  title    = {Nonsymmetric {Reduction}-{Based} {Algebraic} {Multigrid}},
  journal  = {SIAM Journal on Scientific Computing},
  year     = {2019},
  volume   = {41},
  number   = {5},
  pages    = {S242--S268},
  month    = jan,
  issn     = {1064-8275},
  doi      = {10.1137/18M1193761},
  url      = {https://epubs.siam.org/doi/abs/10.1137/18M1193761},
  urldate  = {2020-05-18},
}

@Article{Southworth2017,
  author   = {Manteuffel, Thomas A. and Ruge, John and Southworth, Ben S.},
  title    = {Nonsymmetric Algebraic Multigrid Based on Local Approximate Ideal Restriction ($\ell${AIR})},
  journal  = {SIAM Journal on Scientific Computing},
  year     = {2018},
  volume   = {40},
  number   = {6},
  pages    = {A4105-A4130},
  doi      = {10.1137/17M1144350},
  eprint   = {https://doi.org/10.1137/17M1144350},
  url      = { 
    
        https://doi.org/10.1137/17M1144350
    
    

},
}

@Article{Loe2022,
  author    = {Loe, Jennifer A. and Morgan, Ronald B.},
  title     = {Toward efficient polynomial preconditioning for {GMRES}},
  journal   = {Numerical Linear Algebra with Applications},
  year      = {2022},
  volume    = {29},
  number    = {4},
  pages     = {e2427},
  issn      = {1099-1506},
  copyright = {© 2021 John Wiley \& Sons Ltd.},
  doi       = {10.1002/nla.2427},
  language  = {en},
  url       = {https://onlinelibrary.wiley.com/doi/abs/10.1002/nla.2427},
  urldate   = {2023-08-16},
}

@InProceedings{Luby1985,
  author    = {Luby, M},
  title     = {A simple parallel algorithm for the maximal independent set problem},
  booktitle = {Proceedings of the seventeenth annual {ACM} symposium on {Theory} of computing},
  year      = {1985},
  series    = {{STOC} '85},
  pages     = {1--10},
  address   = {New York, NY, USA},
  month     = dec,
  publisher = {Association for Computing Machinery},
  doi       = {10.1145/22145.22146},
  isbn      = {978-0-89791-151-1},
  url       = {https://dl.acm.org/doi/10.1145/22145.22146},
  urldate   = {2023-11-09},
}

@Article{Heller1976,
  author   = {Heller, Don},
  title    = {Some {Aspects} of the {Cyclic} {Reduction} {Algorithm} for {Block} {Tridiagonal} {Linear} {Systems}},
  journal  = {SIAM Journal on Numerical Analysis},
  year     = {1976},
  volume   = {13},
  number   = {4},
  pages    = {484--496},
  month    = sep,
  issn     = {0036-1429},
  note     = {Publisher: Society for Industrial and Applied Mathematics},
  doi      = {10.1137/0713042},
  url      = {https://epubs.siam.org/doi/abs/10.1137/0713042},
  urldate  = {2023-11-29},
}

@Book{Golub1992,
  title     = {Cyclic reduction/multigrid},
  publisher = {Numerical Analysis Project, Computer Science Department, Stanford University},
  year      = {1992},
  author    = {Golub, Gene H. and Tuminaro, Raymond S.},
  url       = {http://infolab.stanford.edu/pub/cstr/reports/na/m/92/14/NA-M-92-14.pdf},
  urldate   = {2023-11-29},
}

@Article{Dargaville2025,
  author   = {Dargaville, Steven and Smedley-Stevenson, Richard and Smith, Paul and Pain, Christopher C},
  title    = {Coarsening and parallelism with reduction multigrids for hyperbolic {Boltzmann} transport},
  journal  = {The International Journal of High Performance Computing Applications},
  year     = {2025},
  volume   = {39},
  number   = {3},
  pages    = {364--384},
  issn     = {1094-3420},
  doi      = {10.1177/10943420241304759},
  url      = {https://doi.org/10.1177/10943420241304759},
}

@Article{Dargaville2024a,
  author   = {Dargaville, S. and Smedley-Stevenson, R. P. and Smith, P. N. and Pain, C. C.},
  title    = {{AIR} multigrid with {GMRES} polynomials ({AIRG}) and additive preconditioners for {Boltzmann} transport},
  journal  = {Journal of Computational Physics},
  year     = {2024},
  volume   = {518},
  pages    = {113342},
  month    = dec,
  issn     = {0021-9991},
  urldate  = {2024-08-12},
}

@Article{Hanophy2020,
  author   = {Hanophy, Joshua and Southworth, Ben S. and Li, Ruipeng and Manteuffel, Tom and Morel, Jim},
  title    = {Parallel {Approximate} {Ideal} {Restriction} {Multigrid} for {Solving} the {SN} {Transport} {Equations}},
  journal  = {Nuclear Science and Engineering},
  year     = {2020},
  volume   = {194},
  number   = {11},
  pages    = {989--1008},
  month    = jun,
  issn     = {0029-5639},
  urldate  = {2020-09-29},
}

@Article{Manteuffel2019a,
  author   = {Manteuffel, Tom and Southworth, Ben S.},
  title    = {Convergence in {Norm} of {Nonsymmetric} {Algebraic} {Multigrid}},
  journal  = {SIAM Journal on Scientific Computing},
  year     = {2019},
  volume   = {41},
  number   = {5},
  pages    = {S269--S296},
  month    = jan,
  issn     = {1064-8275},
  doi      = {10.1137/18M1193773},
  urldate  = {2020-09-29},
}

@Article{Dargaville2024,
  author   = {Dargaville, S. and Smedley-Stevenson, R. P. and Smith, P. N. and Pain, C. C.},
  title    = {Angular {Adaptivity} in {P0} {Space} and {Reduced} {Tolerance} {Solves} for {Boltzmann} {Transport}},
  journal  = {Nuclear Science and Engineering},
  year     = {2024},
  volume   = {198},
  number   = {6},
  pages    = {1235--1254},
  month    = jun,
  issn     = {0029-5639},
  urldate  = {2024-08-08},
}

@Article{Ali2024,
  author   = {Ali, Ahsan and Brannick, James J. and Kahl, Karsten and Krzysik, Oliver A. and Schroder, Jacob B. and Southworth, Ben S.},
  title    = {Constrained {Local} {Approximate} {Ideal} {Restriction} for {Advection}-{Diffusion} {Problems}},
  journal  = {SIAM Journal on Scientific Computing},
  year     = {2024},
  pages    = {S96--S122},
  month    = may,
  urldate  = {2024-08-02},
}

@Article{Zaman2024,
  author   = {Zaman, Tareq and Nytko, Nicolas and Taghibakhshi, Ali and MacLachlan, Scott and Olson, Luke and West, Matthew},
  title    = {Generalizing reduction-based algebraic multigrid},
  journal  = {Numerical Linear Algebra with Applications},
  year     = {2024},
  volume   = {31},
  number   = {3},
  pages    = {e2543},
  language = {en},
  urldate  = {2024-08-02},
}

@Article{Sivas2021,
  author   = {Sivas, A. A. and Southworth, B. S. and Rhebergen, S.},
  title    = {{AIR} {Algebraic} {Multigrid} for a {Space}-{Time} {Hybridizable} {Discontinuous} {Galerkin} {Discretization} of {Advection}(-{Diffusion})},
  journal  = {SIAM Journal on Scientific Computing},
  year     = {2021},
  volume   = {43},
  number   = {5},
  pages    = {A3393--A3416},
  month    = jan,
  urldate  = {2024-11-25},
}

@Article{Rosa2010,
  author   = {Rosa, Massimiliano and Warsa, James S. and Chang, Jae H.},
  title    = {Fourier {Analysis} of {Inexact} {Parallel} {Block}-{Jacobi} {Splitting} with {Transport} {Synthetic} {Acceleration}},
  journal  = {Nuclear Science and Engineering},
  year     = {2010},
  volume   = {164},
  number   = {3},
  pages    = {248--263},
  month    = mar,
  issn     = {0029-5639},
  doi      = {10.13182/NSE09-26},
  url      = {https://doi.org/10.13182/NSE09-26},
  urldate  = {2023-03-16},
}

@Article{Yavuz1992,
  author   = {Yavuz, Musa and Larsen, Edward W.},
  title    = {Iterative {Methods} for {Solving} x-y {Geometry} {SN} {Problems} on {Parallel} {Architecture} {Computers}},
  journal  = {Nuclear Science and Engineering},
  year     = {1992},
  volume   = {112},
  number   = {1},
  pages    = {32--42},
  month    = sep,
  issn     = {0029-5639},
  doi      = {10.13182/NSE92-A23949},
  url      = {https://doi.org/10.13182/NSE92-A23949},
  urldate  = {2023-03-16},
}

@Article{Adams2020,
  author   = {Adams, Michael P. and Adams, Marvin L. and Hawkins, W. Daryl and Smith, Timmie and Rauchwerger, Lawrence and Amato, Nancy M. and Bailey, Teresa S. and Falgout, Robert D. and Kunen, Adam and Brown, Peter},
  title    = {Provably optimal parallel transport sweeps on semi-structured grids},
  journal  = {Journal of Computational Physics},
  year     = {2020},
  volume   = {407},
  pages    = {109234},
  month    = apr,
  issn     = {0021-9991},
  doi      = {10.1016/j.jcp.2020.109234},
  language = {en},
  url      = {https://www.sciencedirect.com/science/article/pii/S0021999120300085},
  urldate  = {2023-03-16},
}

@InProceedings{Deakin2016,
  author    = {Deakin, Tom and McIntosh-Smith, Simon and Gaudin, Wayne},
  title     = {Many-{Core} {Acceleration} of a {Discrete} {Ordinates} {Transport} {Mini}-{App} at {Extreme} {Scale}},
  booktitle = {High {Performance} {Computing}},
  year      = {2016},
  editor    = {Kunkel, Julian M. and Balaji, Pavan and Dongarra, Jack},
  series    = {Lecture {Notes} in {Computer} {Science}},
  pages     = {429--448},
  address   = {Cham},
  publisher = {Springer International Publishing},
  doi       = {10.1007/978-3-319-41321-1_22},
  isbn      = {978-3-319-41321-1},
  language  = {en},
}

@TechReport{Kunen2019,
  author      = {Kunen, A. and Loffeld, J. and Black, A. and Chen, R. and Nowak, P. and Haut, T. and Bailey, T. and Brown, P. and Rennich, S. and Maginot, P.},
  title       = {Porting {3D} discrete ordinates sweep algorithm in ardra to {CUDA}},
  institution = {Lawrence Livermore National Lab.(LLNL), Livermore, CA (United States)},
  year        = {2019},
  url         = {https://www.osti.gov/servlets/purl/1559411},
  urldate     = {2025-03-10},
}

@Article{Vermaak2021,
  author   = {Vermaak, Jan I. C. and Ragusa, Jean C. and Adams, Marvin L. and Morel, Jim E.},
  title    = {Massively parallel transport sweeps on meshes with cyclic dependencies},
  journal  = {Journal of Computational Physics},
  year     = {2021},
  volume   = {425},
  pages    = {109892},
  month    = jan,
  issn     = {0021-9991},
  doi      = {10.1016/j.jcp.2020.109892},
  url      = {https://www.sciencedirect.com/science/article/pii/S0021999120306665},
  urldate  = {2023-09-13},
}

@Article{Pautz2017,
  author   = {Pautz, Shawn D. and Bailey, Teresa S.},
  title    = {Parallel {Deterministic} {Transport} {Sweeps} of {Structured} and {Unstructured} {Meshes} with {Overloaded} {Mesh} {Decompositions}},
  journal  = {Nuclear Science and Engineering},
  year     = {2017},
  volume   = {185},
  number   = {1},
  pages    = {70--77},
  month    = jan,
  doi      = {10.13182/NSE16-34},
  url      = {https://doi.org/10.13182/NSE16-34},
  urldate  = {2023-09-13},
}

@Article{Balogh2022,
  author    = {Balogh, Gabor D. and Flynn, Tobias S. and Laizet, Sylvain and Mudalige, Gihan R. and Reguly, Istan Z.},
  title     = {Scalable {Many}-{Core} {Algorithms} for {Tridiagonal} {Solvers}},
  journal   = {Computing in Science \& Engineering},
  year      = {2022},
  volume    = {24},
  number    = {1},
  pages     = {26--35},
  month     = jan,
  issn      = {1521-9615, 1558-366X},
  copyright = {https://ieeexplore.ieee.org/Xplorehelp/downloads/license-information/IEEE.html},
  doi       = {10.1109/MCSE.2021.3130544},
  language  = {en},
  url       = {https://ieeexplore.ieee.org/document/9626460/},
  urldate   = {2024-11-27},
}

@Article{Bini2009,
  author   = {Bini, Dario and Meini, Beatrice},
  title    = {The cyclic reduction algorithm},
  journal  = {Numerical Algorithms},
  year     = {2009},
  volume   = {51},
  pages    = {23--60},
  month    = may,
}

@Article{Tolmachev2025,
  author   = {Tolmachev, Dmitrii and Marti, Philippe and Castiglioni, Giacomo and Jackson, Andrew},
  title    = {High {Performance} {Solution} of {Tridiagonal} {Systems} on the {GPU}},
  journal  = {ACM Trans. Parallel Comput.},
  year     = {2025},
  volume   = {12},
  number   = {2},
  pages    = {5:1--5:25},
  month    = may,
  urldate  = {2025-06-27},
}

@Article{Song2022,
  author        = {Song, Hang and Matsuno, Kristen V. and West, Jacob R. and Subramaniam, Akshay and Ghate, Aditya S. and Lele, Sanjiva K.},
  title         = {Scalable parallel linear solver for compact banded systems on heterogeneous architectures},
  journal       = {Journal of Computational Physics},
  year          = {2022},
  volume        = {468},
  pages         = {111443},
  month         = nov,
  __markedentry = {[sdargavi:]},
  abstract      = {A scalable algorithm for solving compact banded linear systems on distributed memory architectures is presented. The proposed method factorizes the original system into two levels of memory hierarchies, and solves it using parallel cyclic reduction on both distributed and shared memory. This method has a lower communication footprint across distributed memory partitions compared to conventional algorithms involving data transposes or re-partitioning. The algorithm developed in this work is generalized to cyclic compact banded systems with flexible data decompositions. For cyclic compact banded systems, the method is a direct solver with a deterministic operation and communication counts depending on the matrix size, its bandwidth, and the partition strategy. The implementation and runtime configuration details are discussed for performance optimization. Scalability is demonstrated on the linear solver as well as on a representative fluid mechanics application problem, in which the dominant computational cost is solving the cyclic tridiagonal linear systems of compact numerical schemes on a 3D periodic domain. The algorithm is particularly useful for solving the linear systems arising from the application of compact finite difference operators to a wide range of partial differential equation problems, such as but not limited to the numerical simulations of compressible turbulent flows, aeroacoustics, elastic–plastic wave propagation, and electromagnetics. It alleviates obstacles to their use on modern high performance computing hardware, where memory and computational power are distributed across nodes with multi-threaded processing units.},
  file          = {ScienceDirect Snapshot:/home/sdargavi/.mozilla/firefox/rmqz5rdm.default/zotero/storage/JXB7S2KL/S0021999122005058.html:text/html;Submitted Version:/home/sdargavi/.mozilla/firefox/rmqz5rdm.default/zotero/storage/NS2GMXIQ/Song et al. - 2022 - Scalable parallel linear solver for compact banded.pdf:application/pdf},
  keywords      = {Compact banded system, Distributed memory, Parallel computing, Parallel cyclic reduction, Periodic boundary},
}

@InProceedings{Gander1998,
  author        = {Gander, Walter and Golub, Gene H.},
  title         = {Cyclic reduction-history and applications},
  booktitle     = {Scientific {Computing}: {Proceedings} of the {Workshop}, 10-12 {March} 1997, {Hong} {Kong}},
  year          = {1998},
  pages         = {73},
  publisher     = {Springer Science \& Business Media},
  __markedentry = {[sdargavi:6]},
  file          = {Available Version (via Google Scholar):/home/sdargavi/.mozilla/firefox/rmqz5rdm.default/zotero/storage/TVEUVC2K/Gander^1 and Golub - 1998 - Cyclic reduction-history and applications.pdf:application/pdf},
  url           = {https://books.google.co.uk/books?hl=en&lr=&id=hzPtYrguXOsC&oi=fnd&pg=PA73&dq=Cyclic+reduction%E2%80%94history+and+applications&ots=gKDyKbUSbR&sig=nCkTNoaZ_K9gxWqDxCor0VX4Nzo},
  urldate       = {2025-07-29},
}

@Misc{petsc-web-page,
  author       = {Satish Balay and Shrirang Abhyankar and Mark~F. Adams and Steven Benson and Jed Brown and Peter Brune and Kris Buschelman and Emil~M. Constantinescu and Lisandro Dalcin and Alp Dener and Victor Eijkhout and Jacob Faibussowitsch and William~D. Gropp and V\'{a}clav Hapla and Tobin Isaac and Pierre Jolivet and Dmitry Karpeev and Dinesh Kaushik and Matthew~G. Knepley and Fande Kong and Scott Kruger and Dave~A. May and Lois Curfman McInnes and Richard Tran Mills and Lawrence Mitchell and Todd Munson and Jose~E. Roman and Karl Rupp and Patrick Sanan and Jason Sarich and Barry~F. Smith and Stefano Zampini and Hong Zhang and Junchao Zhang},
  title        = {{PETS}c {W}eb page},
  howpublished = {\url{https://petsc.org/}},
  year         = {2025},
  url          = {https://petsc.org/},
}

@Article{mills2021,
  author  = {Richard Tran Mills and Mark F. Adams and Satish Balay and Jed Brown and Alp Dener},
  title   = {Toward performance-portable {PETS}c for {GPU}-based exascale systems},
  journal = {Parallel Computing},
  year    = {2021},
  volume  = {108},
  pages   = {102831},
  issn    = {0167-8191},
  doi     = {https://doi.org/10.1016/j.parco.2021.102831},
  url     = {https://www.sciencedirect.com/science/article/pii/S016781912100079X},
}

@Article{pyamg2023,
  author    = {Nathan Bell and Luke N. Olson and Jacob Schroder and Ben Southworth},
  title     = {{PyAMG}: Algebraic Multigrid Solvers in Python},
  journal   = {Journal of Open Source Software},
  year      = {2023},
  volume    = {8},
  number    = {87},
  pages     = {5495},
  doi       = {10.21105/joss.05495},
  publisher = {The Open Journal},
  url       = {https://doi.org/10.21105/joss.05495},
}






% \end{thebibliography}

\end{document}